\documentclass[final,3p,ecuild]{elsarticle}

\usepackage{amssymb}
\usepackage{arydshln}
\usepackage{xcolor}
\usepackage{booktabs}
\usepackage[ruled,vlined,linesnumbered]{algorithm2e}
\usepackage{todonotes}
\usepackage{amsthm}
\usepackage{amsmath}
\usepackage[colorlinks,linkcolor=red,anchorcolor=green,citecolor=blue]{hyperref}
\usepackage{orcidlink}
\usepackage{lineno}
\usepackage{bm}

\begin{document}
\begin{frontmatter}



\title{Function-on-function nonlinear autoregressive model (F2NARX) for probabilistic surrogate modeling of stochastic dynamical systems} 

\author[inst1]{Zhouzhou Song\orcidlink{0000-0002-5899-6205}\corref{cor1}} 
\ead{zhouzhou.song@tu-dortmund.de}

\author[inst1]{Marcos A. Valdebenito\orcidlink{0000-0002-5083-0454}}
\ead{marcos.valdebenito@tu-dortmund.de}

\author[inst2]{Styfen Schär\orcidlink{0000-0001-6715-220X}}
\ead{styfen.schaer@ibk.baug.ethz.ch}

\author[inst2]{Stefano Marelli\orcidlink{0000-0002-9268-9014}}
\ead{marelli@ibk.baug.ethz.ch}

\author[inst2]{Bruno Sudret\orcidlink{0000-0002-9501-7395}}
\ead{sudret@ethz.ch}

\author[inst1,inst3]{Matthias G.R. Faes\orcidlink{0000-0003-3341-3410}}
\ead{matthias.faes@tu-dortmund.de}

\affiliation[inst1]{
            organization={Chair for Reliability Engineering, TU Dortmund University},
            addressline={Leonhard-Euler-Strasse 5},
            postcode={44227 Dortmund}, 
            country={Germany}}

\affiliation[inst2]{
            organization={Chair of Risk, Safety and Uncertainty Quantification, ETH Zurich},
            addressline={Stefano-Franscini-Platz 5},
            postcode={8093, Zürich}, 
            country={Switzerland}}

\affiliation[inst3]{
            organization={International Joint Research Center for Engineering Reliability and Stochastic Mechanics, Tongji University},
            postcode={Shanghai 200092}, 
            country={PR China}}

\cortext[cor1]{Corresponding author.}

\begin{abstract}
Constructing accurate and computationally efficient surrogate models for predicting dynamical system responses is critical in many engineering domains, yet remains challenging due to the strongly nonlinear and high-dimensional mapping from external excitations and system parameters to system responses. 
This work introduces a novel Function-on-Function Nonlinear AutoRegressive model with eXogenous inputs (F2NARX), which reformulates the 
recently proposed $\mathcal{F}$-NARX method from a function-on-function regression perspective.
The proposed framework substantially improves predictive efficiency while maintaining high accuracy.
By combining principal component analysis with Gaussian process regression, F2NARX further enables probabilistic predictions of dynamical responses via the unscented transform in an autoregressive manner.
The effectiveness of the method is demonstrated through case studies of varying complexity.
Results show that F2NARX outperforms state-of-the-art NARX model by orders of magnitude in efficiency while achieving higher accuracy in general.
The probabilistic prediction capability of F2NARX is further demonstrated by combining it with active learning for first-passage probability estimation, showing its potential for reliability-oriented applications.

\end{abstract}




\end{frontmatter}

\section{Introduction}
Dynamical systems are commonly encountered in real-world applications, where they evolve continuously under external excitations, leading to complex time-dependent responses. Accurately and efficiently emulating dynamical systems has attracted increasing attention in various engineering domains, including uncertainty quantification \cite{mai2016surrogate,song2025efficient,cheng2025state}, reliability analysis \cite{li2022lstm,zhou2022efficient,zhang2024rlstm}, design optimization \cite{deshmukh2017design,yang2025structural}, system control \cite{levin1996control,hu2024efficient,chen2025real}, prognostics and health management \cite{theissler2021predictive,de2021framework,feng2023digital}, and digital twins \cite{karkaria2024towards,gunasegaram2024machine}.

In these contexts, constructing accurate surrogate models (sometimes also called emulators) as inexpensive approximations of dynamical systems remains a challenging task. The difficulty arises because the mapping from external excitations and system parameters to system responses is often highly nonlinear, owing to the presence of strongly coupled subsystems or nonlinear components \cite{schar2026mnarx+}. In addition, external excitations are typically high-dimensional time series, and learning such mappings directly suffers from the curse of dimensionality when using common surrogate modeling techniques such as Gaussian processes (or Kriging) \cite{rasmussen2003gaussian,song2024improved}, support vector regression \cite{drucker1996support,roy2023support}, polynomial chaos expansion \cite{xiu2002wiener,Luethen2021Review}, or neural networks \cite{bishop1995neural,chojaczyk2015review}. 
Recent works \cite{simpson2021machine,zhou2022efficient,song2024dimension,zhang2024rlstm,song2025efficient,kim2025dimensionality} employ dimensionality reduction techniques to identify a low-dimensional active (or latent) subspace from high-dimensional inputs and then construct surrogate models within that space. However, effectively identifying the active subspace still requires a large amount of training time histories.

An alternative approach to reducing surrogate modeling complexity is provided by nonlinear autoregressive models with exogenous inputs (NARX) \cite{leontaritis1985input1,leontaritis1985input2,chen1989orthogonal,billings2013nonlinear,mai2016surrogate}. Rather than directly learning the complex mapping from the original excitation to the system response, NARX models assume that the response at a given time instant depends only on its past values, the past and current excitation values, and the system parameters. This local mapping is considerably easier to learn than the original strongly nonlinear and high-dimensional relationship. Another advantage of the NARX framework is that they can incorporate existing surrogate modeling techniques seamlessly to learn the autoregressive mapping. Despite their success in various engineering applications, NARX models face a critical challenge: they require the selection of appropriate time lags for both exogenous inputs and system responses, for which no universal solution is currently available \cite{schar2025surrogate}. Moreover, NARX models may struggle to emulate the responses of strongly nonlinear dynamical systems \cite{cheng2025state}, and when a large number of time lags are required, NARX models remain susceptible to the curse of dimensionality.

In response to this limitation,  a manifold-based variant of NARX, namely manifold-NARX (mNARX) \cite{schar2024emulating} was recently introduced to emulate highly complex dynamical systems. In this approach, a problem-specific exogenous input manifold is sequentially constructed from a combination of available input and output quantities, to serve as a more suitable basis for constructing a multi-dimensional NARX model. However, constructing the manifold requires expert knowledge on the problem physics, as well as the expected relations between different model responses. 
Despite its increased expressivity, manifold-based methodology still borrows its core learning mechanism from classical NARX, thus exhibiting a similarly strong dependence both on the choice of input lags, and on the time-discretization characteristics of the problem.
To address this issue, a novel approach to NARX, namely functional-NARX ($\mathcal{F}$-NARX) \cite{schar2025surrogate}, was recently developed.
Instead of directly relying on the values of exogenous inputs and responses at discretized time lags, $\mathcal{F}$-NARX employs time-dependent functional features (although possibly discretized) extracted from a local time window, which more efficiently represents the dynamics of both the model inputs and its responses. 
This approach almost entirely eliminates the sensitivity of NARX to time-discretization, as well as strongly reducing its well known \textit{over-reliance problem} \cite{piroddi2003identification}.
Recent work also combines both mNARX and $\mathcal{F}$-NARX to drastically reduce the amount of expert knowledge required, by enabling a semi-automatic manifold identification strategy \cite{schar2026mnarx+}.

However, none of these methods offers an efficient way of producing probabilistic predictions, such as confidence bounds or uncertainty envelopes, which are crucial for robust estimation and decision-making. 
In addition, such advanced NARX models still face computational-efficiency challenges when long-term predictions or many surrogate-model runs are required, because they rely on a one-step-ahead prediction strategy.

In this work, we propose a novel \textit{function-on-function nonlinear autoregressive model with exogenous inputs} (F2NARX) for efficient emulation of complex dynamical systems.
The proposed model can be viewed as an extension of $\mathcal{F}$-NARX in terms of the autoregressive structure.
Instead of following the classical one-step-ahead prediction of NARX models, F2NARX adopts a one-time-window-ahead strategy. 
It assumes that the response within a given time window depends on the response from the previous window, the excitation values in both the considered and previous windows, and the system parameters, thereby shifting towards a function-on-function mapping. 
We adopt the same principal component analysis (PCA)-based discretization approach as in $\mathcal{F}$-NARX to extract features from excitation and response functions within local time windows. 
The original function-on-function mapping is then discretized into a multi-input multi-output (MIMO) mapping, which is further decomposed into a set of single-output mappings by exploiting the orthogonality of PCA.
Sparse Gaussian process (SGP) regression is used to learn these single-output mappings from large autoregressive training data sets.
Meanwhile, we develop a probabilistic prediction scheme that combines the predictive uncertainty of SGPs with the unscented transform to quantify the epistemic prediction uncertainty of dynamical responses in an autoregressive manner.

While F2NARX can be applied to a wide range of applications, this work focuses on probabilistic surrogate modeling of stochastic dynamical systems.
Case studies of varying complexity demonstrate that F2NARX can emulate dynamical responses with up to orders-of-magnitude reductions in computational time compared with state-of-the-art NARX models, while maintaining high accuracy. 
The usefulness of the probabilistic predictions is further demonstrated through active learning-based first-passage probability estimation.

The remainder of this paper is organized as follows. 
Section \ref{sec_2} presents background information on nonlinear dynamical systems and the NARX model. 
Section \ref{sec_3} introduces the proposed F2NARX method, including its model formulation, training procedure, and probabilistic prediction scheme. 
Section \ref{sec_AL} demonstrates the use of probabilistic F2NARX predictions in an active learning-based first-passage probability estimation problem.
Section \ref{sec_4} reports two case studies of varying complexity to demonstrate the effectiveness of F2NARX.
Finally, Section \ref{sec_5} discusses the main findings of the paper and provides concluding remarks.

\section{Background}
\label{sec_2}
This work focuses on applying NARX-based methods to probabilistic surrogate modeling of stochastic structural dynamical systems. In this section, we first introduce the formulation of structural dynamical systems with uncertain system and excitation parameters, followed by a brief overview of conventional and advanced NARX models.

\subsection{Stochastic structural dynamical systems}
The motion of an $N_d$-degree-of-freedom nonlinear stochastic dynamical system is modeled by the following equation:
\begin{equation}
    \label{governing_eq_MDOF}
    \bm{M}(\bm{\varTheta})\ddot{\bm{Y}}(t)+\bm{C}(\bm{\varTheta})\dot{\bm{Y}}(t)+\bm{R}(\dot{\bm{Y}}(t),\bm{Y}(t),\bm{\varTheta})=\bm{U}(\bm{\varTheta},\bm{\varPhi},t),
\end{equation}
where $\bm{M}$ and $\bm{C}$ are the $N_d\times N_d$ mass and damping matrices, respectively, both dependent on the parameter vector $\bm{\varTheta}$; $\ddot{\bm{Y}}(t)$, $\dot{\bm{Y}}(t)$, and $\bm{Y}(t)$ are the $N_d$-dimensional acceleration, velocity, and displacement vectors, respectively; $\bm{R}(\dot{\bm{Y}}(t),\bm{Y}(t),\bm{\varTheta})$ is the $N_d$-dimensional nonlinear restoring force vector, which is a nonlinear function of $\dot{\bm{Y}}(t)$, $\bm{Y}(t)$, and $\bm{\varTheta}$; $\bm{U}(\bm{\varTheta},\bm{\varPhi},t)$ denotes the $N_d$-dimensional external excitation vector, which depends on both parameter vector $\bm{\varTheta}$ and the high-dimensional random vector $\bm{\varPhi}$ that controls the randomness of the excitation. Given fixed values $\bm{\theta}$ and $\bm{\phi}$ for $\bm{\varTheta}$ and $\bm{\varPhi}$, respectively, the excitation $\bm{u}(\bm{\theta},\bm{\phi},t)$ is uniquely determined as a time series.
Then, for a given time interval $[t_0,t_e]$ and initial conditions $\dot{\bm{Y}}(t_0)=\dot{\bm{y}}_0,\ \bm{Y}(t_0)=\bm{y}_0$, the dynamic response $\bm{y}(t)$ is uniquely determined over $[t_0,t_e]$ by numerically solving Eq.~(\ref{governing_eq_MDOF}), for example using the Runge-Kutta method.

\subsection{Nonlinear autoregressive models with exogenous inputs}
The NARX model is a widely used method for surrogate modeling of stochastic dynamical systems. Its core idea is to capture the system's dynamics via a discrete-time-step representation, which expresses the system response at a given time instant as a function of its past values, past and current excitation values, and key parameters. Taking a single response quantity of interest (QoI) $y(t)$ as an example, the NARX model is expressed as:
\begin{equation}
    \label{NARX_model}
    y(t^*)=\mathcal{F}(\bm{u}(t^*),\bm{u}(t^*-\delta t),\cdots,\bm{u}(t^*-n_u\delta t),y(t^*-\delta t),y(t^*-2\delta t),\cdots,y(t^*-n_y\delta t),\bm{\varTheta}),
\end{equation}
where $\mathcal{F}(\cdot)$ is the underlying function to be learned, $\delta t$ is a small positive time increment, and $n_u$ and $n_y$ are the maximum excitation and response time lags. The function $\mathcal{F}(\cdot)$ is typically approximated using polynomials or Gaussian processes. However, this explicit discrete-time-step representation has limitations, including the difficulty in determining suitable values for $n_u$ and $n_y$ \cite{cheng2025state}, poor performance in predicting strongly nonlinear dynamical systems \cite{cheng2025state}, and the curse of dimensionality \cite{schar2025surrogate}.

To overcome these limitations, Schär et al. \cite{schar2025surrogate} recently proposed treating the NARX model from a continuous functional perspective, referred to as $\mathcal{F}$-NARX. The $\mathcal{F}$-NARX model represents the response at a given time instant $y(t^*)$ as a function of the response in the previous time window $y(\tau_1-T),\tau_1\in[t^*,t^*+T-\delta t]$, the excitation over a time window including both previous and current values $\bm{u}(\tau_2-T),\tau_2\in[t^*,t^*+T]$, and the parameters $\bm{\varTheta}$.
Following this, the $\mathcal{F}$-NARX model can be expressed as:
\begin{equation}
    \label{F-NARX model}
    y(t^*)=\mathcal{F}(\bm{u}(\tau_2-T),y(\tau_1-T),\bm{\varTheta}),\tau_1\in[t^*,t^*+T-\delta t],\tau_2\in[t^*,t^*+T].
\end{equation}
In practice, $\bm{u}(\tau_2-T)$ and $y(\tau_1-T)$ are discretized into low-dimensional feature vectors obtained via feature-extraction methods such as principal component analysis (PCA). These features represent the dynamics more efficiently and help mitigate the sensitivity of NARX to time discretization.
Although $\mathcal{F}$-NARX has demonstrated advantages in both accuracy and robustness over the classical NARX model, its prediction efficiency remains limited because it relies on a one-step-ahead prediction strategy.
As a result, the computational cost can become relatively large when long-term predictions or many surrogate-model runs are required.
Moreover, an efficient scheme for probabilistic prediction remains unavailable.

\section{Probabilistic function-on-function nonlinear autoregressive model for emulating dynamical systems}
\label{sec_3}
In this section, we present the proposed method for emulating nonlinear dynamical systems. We begin by revisiting the $\mathcal{F}$-NARX model \cite{schar2025surrogate} from its function-on-scalar perspective to a function-on-function one. Next, we introduce the training procedure of the F2NARX method, which includes training data generation, functional feature extraction, and the construction of a sparse Gaussian process regression-based emulator. Subsequently, we describe the procedure for probabilistic prediction of dynamical responses using the F2NARX model. 

\subsection{Function-on-function nonlinear autoregressive model with exogenous inputs (F2NARX)}
\label{subsec3.1}

\begin{figure}[t]
\centering
\includegraphics[scale=0.45]{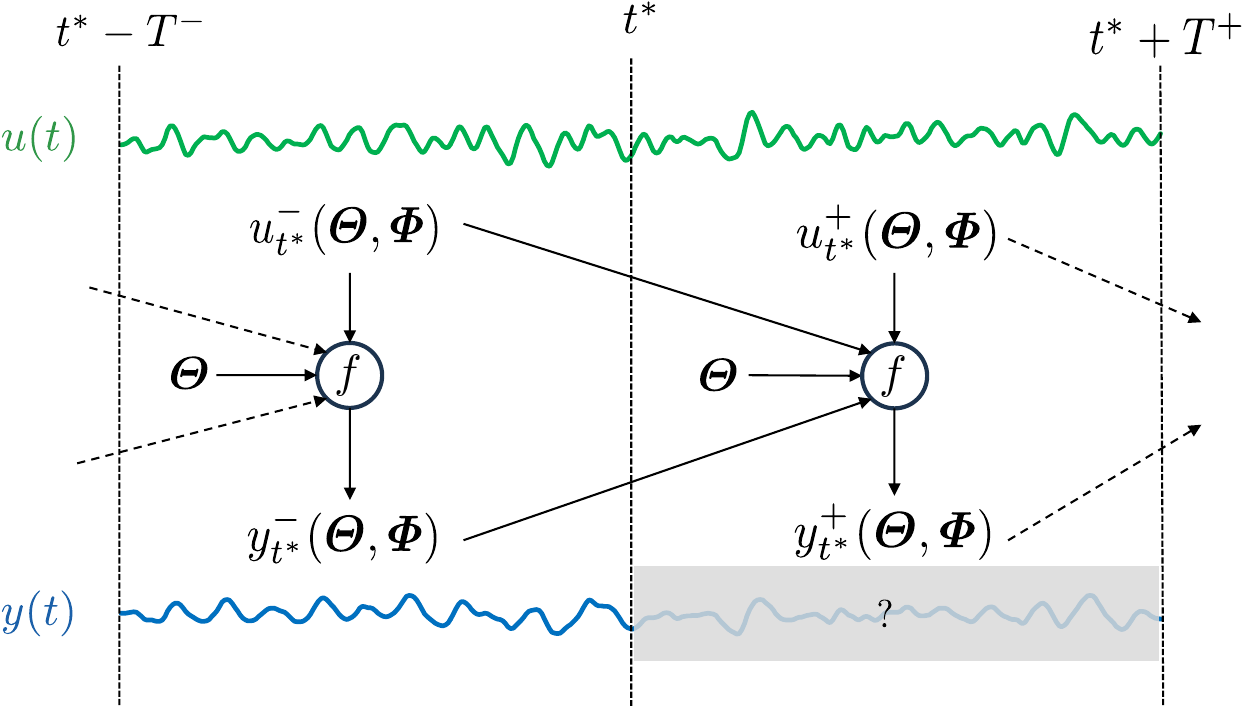}
\caption{The illustration of F2NARX model.}
\label{fig_F2NARX}
\end{figure}

For simplicity, this study focuses on the case with a single external excitation. As illustrated in Fig.~\ref{fig_F2NARX}, the excitation and response functions are segmented into local time windows. 
Given a reference time instant $t^*$, we define the following auxiliary functions:
\begin{equation}
    \label{eqn:aux functions 2}
    \begin{split}
        y^{-}_{t^*}(\bm{\varTheta},\bm{\varPhi}) \equiv&~ y(\bm{\varTheta},\bm{\varPhi},t),~~~~ t \in \left(t^* - T^-, t^* \right],\\
        y^{+}_{t^*}(\bm{\varTheta},\bm{\varPhi}) \equiv&~ y(\bm{\varTheta},\bm{\varPhi},t),~~~~ t \in \left(t^*, t^* + T^+ \right],\\
        u^{-}_{t^*}(\bm{\varTheta},\bm{\varPhi}) \equiv&~ u(\bm{\varTheta},\bm{\varPhi},t),~~~~ t \in \left(t^* - T^-, t^* \right],\\
        u^{+}_{t^*}(\bm{\varTheta},\bm{\varPhi}) \equiv&~ u(\bm{\varTheta},\bm{\varPhi},t),~~~~ t \in \left(t^*, t^* + T^+ \right],
    \end{split}
\end{equation}
where $y^{\mp}_{t^*}(\bm{\varTheta},\bm{\varPhi})$ are the response quantity of interest $y(\bm{\varTheta},\bm{\varPhi},t)$ on a time window of length $T^-$ before $t^*$, and on a time window of length $T^+$ after $t^*$, respectively; $u^{\mp}_{t^*}(\bm{\varTheta},\bm{\varPhi})$ are the excitation $u(\bm{\varTheta},\bm{\varPhi},t)$ on a time window of length $T^-$ before $t^*$, and on a time window of length $T^+$ after $t^*$, respectively. For simplicity, we denote $y^{\mp}_{t^*}(\bm{\varTheta},\bm{\varPhi})$ and $u^{\mp}_{t^*}(\bm{\varTheta},\bm{\varPhi})$ as $y^{\mp}_{t^*}$ and $u^{\mp}_{t^*}$, respectively.
The core idea of the F2NARX model is to represent $y^{+}_{t^*}$ as a function of $u^{+}_{t^*}$, $u^{-}_{t^*}$, $y^{-}_{t^*}$, and $\bm{\varTheta}$. This relationship is formulated as:
\begin{equation}
\label{eq_F2NARX_general}
y^{+}_{t^*}=f(u^{+}_{t^*},u^{-}_{t^*},y^{-}_{t^*},\bm{\varTheta}).
\end{equation}
Here we assume that $y^{+}_{t^*}$, $y^{-}_{t^*}$, $u^{+}_{t^*}$, and $u^{-}_{t^*}$ are square-integrable functions defined on their respective time windows. Therefore, Eq.~(\ref{eq_F2NARX_general}) captures the evolution of the dynamical system via a function-on-function mapping.
Note that the model in Eq.~(\ref{eq_F2NARX_general}) can be readily extended to scenarios with multiple excitation functions by replacing the scalar function $u$ with a vector-valued function $\bm{u}$.

The window width $T^-$, referred to as the memory of the model, defines the \textit{look-back time} \cite{schar2025surrogate}. Similarly, the prediction window width $T^+$ is defined as the \textit{look-ahead time} in this study. These parameters carry physical meaning: the system response in the time window $(t^*-T^-,t^*]$ significantly influences the response in the time window $(t^*,t^*+T^+]$. In general, different look-back times may be assigned to $y(t)$ and $u(t)$, and the look-ahead time does not necessarily need to match the look-back time. However, to simplify the formulation without loss of generality, we assume that the look-ahead time and all look-back times are the same, i.e., $T^+=T^-=T$, throughout the remainder of this paper. 

It should be noted that when $t^*=t_0$, $(t_0,t_0+T]$ corresponds to the first local time window, for which no previous window exists. In this case, the model in Eq.~(\ref{eq_F2NARX_general}) is reduced to:
\begin{equation}
\label{eq_F2NARX_first_window}
y_{t_0}^+ = f_0(u_{t_0}^+, u(t_0), y(t_0), \bm{\varTheta}),
\end{equation}
where $u(t_0)$ and $y(t_0)$ denote the initial values of $u(t)$ and $y(t)$, respectively.
Note that, in Eq.~\eqref{eq_F2NARX_first_window}, the initial values are treated separately from the uncertain parameter vector $\bm{\varTheta}$ and are therefore not included in $\bm{\varTheta}$.

\subsection{F2NARX model training}
\label{subsec3.2}

\subsubsection{Training data set generation}
In practice, learning a function-on-function mapping directly is challenging, as a function is inherently infinite-dimensional and is numerically represented by a high-dimensional vector. 
Alternatively, a more practical approach is to extract a small number of features from the functions and to construct a surrogate model on these features, thereby extending the $\mathcal{F}$-NARX approach to predict the features of the current window. This strategy enables a transition from a continuous function-to-function mapping to a much more manageable feature-to-feature mapping.
We define $\mathcal{K}_{u}$ and $\mathcal{K}_{y}$ as two mappings that extract features from $u^{\mp}_{t^*}$ and $y^{\mp}_{t^*}$, respectively, whose technical details will be presented in Section~\ref{feature extraction}.
Additionally, we denote $\bm{\xi}_{u_{t^*}^+}\in
\mathbb{R}^{1\times m_{u}}$, $\bm{\xi}_{u_{t^*}^-}\in
\mathbb{R}^{1\times m_{u}}$, $\bm{\xi}_{y_{t^*}^+}\in
\mathbb{R}^{1\times m_{y}}$, and $\bm{\xi}_{y_{t^*}^-}\in
\mathbb{R}^{1\times m_{y}}$ as the feature vectors of $u_{t^*}^+$, $u_{t^*}^-$, $y_{t^*}^+$, and $y_{t^*}^-$, where $m_{u}$ and $m_{y}$ represent the numbers of features for the local excitation and response functions, respectively. The feature vectors can be obtained as:
\begin{equation}
\label{eq-feature-extraction-mappings}
\bm{\xi}_{u_{t^*}^+}=\mathcal{K}_{u}(u_{t^*}^+), \quad \bm{\xi}_{u_{t^*}^-}=\mathcal{K}_{u}(u_{t^*}^-), \quad \bm{\xi}_{y_{t^*}^+}=\mathcal{K}_{y}(y_{t^*}^+), \quad \bm{\xi}_{y_{t^*}^-}=\mathcal{K}_{y}(y_{t^*}^-).
\end{equation}
Based on Eq.~(\ref{eq-feature-extraction-mappings}), Eq.~(\ref{eq_F2NARX_general}) can be transformed into an MIMO mapping:
\begin{equation}
\label{eq_MIMO}
\bm{\xi}_{y_{t^*}^+} = \bm{f}(\bm{\xi}_{u_{t^*}^+}, \bm{\xi}_{u_{t^*}^-}, \bm{\xi}_{y_{t^*}^-}, \bm{\varTheta}).
\end{equation}
Similarly, Eq.~(\ref{eq_F2NARX_first_window}) can be transformed into:
\begin{equation}
\label{eq_MIMO_first_window}
\bm{\xi}_{y_{t_0}^+} = \bm{f}_0(\bm{\xi}_{u_{t_0}^+}, u(t_0), y(t_0), \bm{\varTheta}).
\end{equation}

To construct a surrogate model for the system described in Eq.~(\ref{eq_MIMO}), the training data set must be generated. 
First, $N_{\text{ED}}$ samples $\left\{\bm{\theta}^{(i)}\right\}_{i=1}^{N_{\text{ED}}}$ and $\left\{\bm{\phi}^{(i)}\right\}_{i=1}^{N_{\text{ED}}}$ of $\bm{\varTheta}$ and $\bm{\varPhi}$ are generated by using Monte Carlo simulation (MCS). 
Based on this, $N_{\text{ED}}$ realizations of the excitation function $u^{(1)}(t), \dots, u^{(N_{\text{ED}})}(t)$ can be obtained. 
By either running the computational model of the dynamical system or utilizing observed experimental data, $N_{\text{ED}}$ realizations of the response function $y^{(1)}(t), \dots, y^{(N_{\text{ED}})}(t)$ are acquired. 
In practice, $u^{(1)}(t), \dots, u^{(N_{\text{ED}})}(t)$ and $y^{(1)}(t), \dots, y^{(N_{\text{ED}})}(t)$ are discretized into vectors using a time increment $\delta t$. 
Consequently, the trajectories of the excitation and response functions can be represented in terms of the following row vectors:
\begin{equation}
\begin{gathered}
u^{(i)}(t) = \bm{u}^{(i)} = [u^{(i)}(t_0), u^{(i)}(t_0 + \delta t), \cdots, u^{(i)}(t_0 + (N_t - 1)\delta t)], \quad i = 1, \cdots, N_{\text{ED}},\\
y^{(i)}(t) = \bm{y}^{(i)} = [y^{(i)}(t_0), y^{(i)}(t_0 + \delta t), \cdots, y^{(i)}(t_0 + (N_t - 1)\delta t)], \quad i = 1, \cdots, N_{\text{ED}},
\end{gathered}
\end{equation}
where $N_t$ is the total number of time instants for one trajectory. Note that, although all training time histories have the same length in this study, the proposed method is also capable of handling data with varying lengths.

Suppose each time window contains $n_T$ time instants and the entire time period $[t_0, t_e]$ is divided into $n_W$ windows. 
If $N_t-1$ is divisible by $n_T$, the time period can be partitioned into $n_W=(N_t-1)/n_T$ non-overlapping windows. 
Otherwise, if $N_t-1$ is not divisible by $n_T$, $n_W$ is set to $\lfloor (N_t-1)/n_T \rfloor + 1$, where $\lfloor \cdot \rfloor$ denotes the floor operator. In this case, the first $\lfloor (N_t-1)/n_T \rfloor$ windows are non-overlapping, while the last two windows overlap. 
It is worth noting that overlapping time windows could also be employed to construct the training dataset; however, the present work focuses on non-overlapping windows. 
We denote the local excitation and response functions in the $j$-th time window for the $i$-th experimental design as $\tilde{u}^{(i)}_j(t)$ and $\tilde{y}^{(i)}_j(t)$ respectively, and use $\tilde{\bm{u}}^{(i)}_j$ and $\tilde{\bm{y}}^{(i)}_j$ to represent their discretized version.
Then, we can obtain two $(N_\text{ED}\cdot n_W) \times n_T$ matrices, $\tilde{\bm{U}}$ and $\tilde{\bm{Y}}$, containing all the local excitation and response functions:
\begin{equation}
\label{matrix for feature extraction}
\tilde{\bm{U}} = 
\begin{pmatrix}
\tilde{\bm{u}}^{(1)}_1 \\
\vdots \\
\tilde{\bm{u}}^{(1)}_{n_W} \\[0.5em]
\hdashline 
\vdots \\[0.5em]
\hdashline 
\tilde{\bm{u}}^{(N_{\text{ED}})}_1 \\
\vdots \\
\tilde{\bm{u}}^{(N_{\text{ED}})}_{n_W}
\end{pmatrix}, \quad
\tilde{\bm{Y}} = 
\begin{pmatrix}
\tilde{\bm{y}}^{(1)}_1 \\
\vdots \\
\tilde{\bm{y}}^{(1)}_{n_W} \\[0.5em]
\hdashline 
\vdots \\[0.5em]
\hdashline 
\tilde{\bm{y}}^{(N_{\text{ED}})}_{1} \\
\vdots \\
\tilde{\bm{y}}^{(N_{\text{ED}})}_{n_W}
\end{pmatrix}.
\end{equation}

After performing feature extraction (discussed in detail in Section~\ref{feature extraction}), $\tilde{\bm{U}}$ is reduced to an ${(N_\text{ED}\cdot n_W) \times m_{u}}$ matrix $\bm{\varXi}_{u}$, and $\tilde{\bm{Y}}$ is reduced to an $(N_\text{ED}\cdot n_W) \times m_{y}$ matrix $\bm{\varXi}_{y}$:
\begin{equation}
\bm{\varXi}_{u} =
\begin{pmatrix}
\mathcal{K}_{u}(\tilde{\bm{u}}^{(1)}_1) \\
\vdots \\
\mathcal{K}_{u}(\tilde{\bm{u}}^{(1)}_{n_W}) \\[0.5em]
\hdashline 
\vdots \\[0.5em]
\hdashline 
\mathcal{K}_{u}(\tilde{\bm{u}}^{(N_{\text{ED}})}_1) \\
\vdots \\
\mathcal{K}_{u}(\tilde{\bm{u}}^{(N_{\text{ED}})}_{n_W})
\end{pmatrix}=
\begin{pmatrix}
\bm{\xi}^{(1)}_{u,1} \\
\vdots \\
\bm{\xi}^{(1)}_{u,n_W} \\[0.5em]
\hdashline 
\vdots \\[0.5em]
\hdashline 
\bm{\xi}^{(N_{\text{ED}})}_{u,1} \\
\vdots \\
\bm{\xi}^{(N_{\text{ED}})}_{u,n_W}
\end{pmatrix}, \quad
\bm{\varXi}_{y}=
\begin{pmatrix}
\mathcal{K}_{y}(\tilde{\bm{y}}^{(1)}_1) \\
\vdots \\
\mathcal{K}_{y}(\tilde{\bm{y}}^{(1)}_{n_W}) \\[0.5em]
\hdashline 
\vdots \\[0.5em]
\hdashline 
\mathcal{K}_{y}(\tilde{\bm{y}}^{(N_{\text{ED}})}_1) \\
\vdots \\
\mathcal{K}_{y}(\tilde{\bm{y}}^{(N_{\text{ED}})}_{n_W})
\end{pmatrix}=
\begin{pmatrix}
\bm{\xi}^{(1)}_{y,1} \\
\vdots \\
\bm{\xi}^{(1)}_{y,n_W} \\[0.5em]
\hdashline 
\vdots \\[0.5em]
\hdashline 
\bm{\xi}^{(N_{\text{ED}})}_{y,1} \\
\vdots \\
\bm{\xi}^{(N_{\text{ED}})}_{y,n_W}
\end{pmatrix},
\end{equation}
in which $\bm{\xi}^{(i)}_{u,j}$ and $\bm{\xi}^{(i)}_{y,j}$ represent the feature vectors of the excitation and response vectors in the $j$-th time window for the $i$-th experimental design, respectively.
Note that, in the current work, feature extraction is performed over the entire collection of windowed excitations or responses.
This treatment relies on the assumption that the windowed excitations or responses can be well approximated by a shared low-dimensional subspace.
However, this assumption may not hold for strongly non-stationary systems, either in amplitude or frequency, whose dominant excitation or response patterns change substantially over time.
In this regard, window-specific feature extraction or more advanced feature extraction methods capable of capturing strong non-stationarity may be considered as possible alternatives.

Consequently, the input and output data used to train the mapping $\bm{f}_0(\cdot)$ are as follows:
\begin{equation}
\label{training_data_f_0}
\bm{X}_0^{\text{train}} =
\begin{pmatrix}
\bm{\xi}^{(1)}_{u,1} & u^{(1)}(t_0) & y^{(1)}(t_0) & \bm{\theta}^{(1)} \\
\vdots & \vdots & \vdots & \vdots \\
\bm{\xi}^{(N_{\text{ED}})}_{u,1} & u^{(N_{\text{ED}})}(t_0) & y^{(N_{\text{ED}})}(t_0) & \bm{\theta}^{(N_{\text{ED}})}
\end{pmatrix}, \quad
\bm{\varXi}_{y,0}^{\text{train}} =
\begin{pmatrix}
\bm{\xi}^{(1)}_{y,1} \\
\vdots \\
\bm{\xi}^{(N_{\text{ED}})}_{y,1}
\end{pmatrix},
\end{equation}
where $\bm{\theta}^{(i)}, i = 1, \dots, N_\text{ED}$ are the experimental designs for the system uncertain parameters. 
This provides $N_\text{ED}$ pairs of input and output data for training $\bm{f}_0(\cdot)$. 
The input and output data used to train the mapping $\bm{f}(\cdot)$ are as follows:
\begin{equation}
\label{training_data_f}
\bm{X}^{\text{train}} = 
\begin{pmatrix}
\bm{\xi}^{(1)}_{u,2} & \bm{\xi}^{(1)}_{u,1} & \bm{\xi}^{(1)}_{y,1} & \bm{\theta}^{(1)} \\
\vdots & \vdots & \vdots & \vdots \\
\bm{\xi}^{(1)}_{u,n_W} & \bm{\xi}^{(1)}_{u,n_W-1} & \bm{\xi}^{(1)}_{y,n_W-1} & \bm{\theta}^{(1)} \\[0.5em]
\hdashline 
\vdots & \vdots & \vdots & \vdots \\[0.5em]
\hdashline 
\bm{\xi}^{(N_{\text{ED}})}_{u,2} & \bm{\xi}^{(N_{\text{ED}})}_{u,1} & \bm{\xi}^{(N_{\text{ED}})}_{y,1} & \bm{\theta}^{(N_{\text{ED}})} \\
\vdots & \vdots & \vdots & \vdots \\
\bm{\xi}^{(N_{\text{ED}})}_{u,n_W} & \bm{\xi}^{(N_{\text{ED}})}_{u,n_W-1} & \bm{\xi}^{(N_{\text{ED}})}_{y,n_W-1} & \bm{\theta}^{(N_{\text{ED}})}
\end{pmatrix}, \quad
\bm{\varXi}_{y}^{\text{train}} = 
\begin{pmatrix}
\bm{\xi}^{(1)}_{y,2} \\
\vdots \\
\bm{\xi}^{(1)}_{y,n_W} \\[0.5em]
\hdashline 
\vdots \\[0.5em]
\hdashline 
\bm{\xi}^{(N_{\text{ED}})}_{y,2} \\
\vdots \\
\bm{\xi}^{(N_{\text{ED}})}_{y,n_W}
\end{pmatrix}.
\end{equation}
This provides $N_\text{ED} \cdot(n_W - 1)$ pairs of input and output data for training $\bm{f}(\cdot)$.

\subsubsection{Feature extraction of local excitation and response functions}
\label{feature extraction}
Many feature extraction or dimensionality reduction methods can be used to obtain the feature mappings $\mathcal{K}_{u}$ and $\mathcal{K}_{y}$, such as principal component analysis (PCA) \cite{pearson1901liii}, independent component analysis \cite{hyvarinen2009independent}, autoencoders \cite{rumelhart1986parallel}, among others. In this research, we employ PCA for feature extraction due to its simplicity and ease of implementation.

In the remainder of this section, we use the local response function $y(t)$ as an example to illustrate the feature extraction process, as the local excitation function follows the same procedure. Let us consider the $(N_\text{ED}\cdot n_W)\times n_T$ matrix $\tilde{\bm{Y}}$, which contains the discretized values of the local response functions, as shown in Eq.~(\ref{matrix for feature extraction}).
The goal of PCA is to find an $n_T\times m_{y}$ orthonormal projection matrix $\bm{V}_{y}$ that extracts low-dimensional features $\bm{\varXi}_{y}$ from $\tilde{\bm{Y}}$ as:
\begin{equation}
    \bm{\varXi}_{y}=\mathcal{K}_{y}(\tilde{\bm{Y}})=\tilde{\bm{Y}}\bm{V}_{y}.
\end{equation}
To obtain $\bm{V}_{y}$, first compute:
\begin{equation}
    \bm{C}=\frac{1}{N_{\text{ED}}\cdot n_W-1}\bm{\tilde{\bm{Y}}}^\top\bm{\tilde{\bm{Y}}}.
\end{equation}
Note that sometimes standardization of $\tilde{\bm{Y}}$ is required.
Subsequently, the eigenvalue decomposition is performed on $\bm{C}$ to obtain the eigenvalues $\lambda_1\geq\lambda_2\geq\cdots\geq \lambda_{n_T}$ and their corresponding eigenvectors $\bm{v}_1,\bm{v}_2,\cdots,\bm{v}_{n_T}$.

In practice, only the first few eigenvectors ${\bm{v}_1,\bm{v}_2,\cdots,\bm{v}_{m_{y}}}$ are sufficient to represent $\bm{C}$. Here, the variance proportion-based criterion is utilized to determine $m_{y}$. That is, $m_{y}$ is chosen as the smallest value that satisfies:
\begin{equation}
    \frac{\sum_{i=1}^{m_{y}} \lambda_i}{\sum_{i=1}^{n_T} \lambda_i}\geq \varepsilon_{\lambda},
\end{equation}
where $\varepsilon_{\lambda}$ is a user-defined threshold, e.g., 0.9999. The feature extraction matrix $\bm{V}_{y}$ is formed by collecting the first $m_{y}$ eigenvectors as $\bm{V}_{y}=[\bm{v}_1,\bm{v}_2,\cdots,\bm{v}_{m_{y}}]$. Accordingly, the feature mapping $\mathcal{K}_{y}$ is defined as:
\begin{equation}
\label{feature-mapping-expression}
\bm{\xi}^{(i)}_{y,j}=\mathcal{K}_{y}\left(\tilde{\bm{y}}^{(i)}_{j}\right)=\tilde{\bm{y}}^{(i)}_{j}\bm{V}_{y}.
\end{equation}
The corresponding inverse mapping $\mathcal{K}_{y}^{-1}$ is given by:
\begin{equation}
\label{inverse-feature-mapping-expression}
    \tilde{\bm{y}}^{(i)}_{j}=\mathcal{K}_{y}^{-1}(\bm{\xi}^{(i)}_{y,j})=\bm{\xi}^{(i)}_{y,j}\bm{V}_{y}^\top.
\end{equation}

\subsubsection{Gaussian process-based emulator for learning F2NARX model}
If we construct a separate surrogate model for each component of $\bm{\xi}_{y_{t^*}^+}$, the two MIMO mappings in Eq.~(\ref{eq_MIMO}) and Eq.~(\ref{eq_MIMO_first_window}) can be further written as two series of single-output mappings respectively as:
\begin{equation}
\label{eq_MISO}
\begin{cases}
        \xi_{y_{t^*}^+}^{(1)} = f^{(1)}(\bm{\xi}_{u_{t^*}^+}, \bm{\xi}_{u_{t^*}^-}, \bm{\xi}_{y_{t^*}^-}, \bm{\varTheta}), \\
        \xi_{y_{t^*}^+}^{(2)} = f^{(2)}(\bm{\xi}_{u_{t^*}^+}, \bm{\xi}_{u_{t^*}^-}, \bm{\xi}_{y_{t^*}^-}, \bm{\varTheta}), \\
        \vdots \\
        \xi_{y_{t^*}^+}^{(m_y)} = f^{(m_y)}(\bm{\xi}_{u_{t^*}^+}, \bm{\xi}_{u_{t^*}^-}, \bm{\xi}_{y_{t^*}^-}, \bm{\varTheta}),
\end{cases}
\end{equation}

\begin{equation}
\label{eq_MISO_first_window}
\begin{cases}
        \xi_{y_{t_0}^+}^{(1)} = f^{(1)}(\bm{\xi}_{u_{t_0}^+}, u(t_0), y(t_0), \bm{\varTheta}), \\
        \xi_{y_{t_0}^+}^{(2)} = f^{(2)}(\bm{\xi}_{u_{t_0}^+}, u(t_0), y(t_0), \bm{\varTheta}), \\
        \vdots \\
        \xi_{y_{t_0}^+}^{(m_y)} = f^{(m_y)}(\bm{\xi}_{u_{t_0}^+}, u(t_0), y(t_0), \bm{\varTheta}).
\end{cases}
\end{equation}

We utilize a Gaussian process (GP) regression model \cite{rasmussen2003gaussian} to emulate each quantity of $\bm{\xi}_{y_{t^*}^+}$. For simplicity, we denote $\xi$ as one of the quantities of $\bm{\xi}_{y_{t^*}^+}$ and denote $\bm{x}\in\mathbb{R}^{1\times p}$ as the model inputs, where $p=n_s+m_{u}+2$ for $f_0$ and $p=n_s+2m_{u}+m_{y}$ for $f$.
A trained GP can provide the predictive mean and variance of $\xi$ at any $\bm{x}^*$ as $\mu_{\hat{\xi}}(\bm{x}^*)$ and $\sigma_{\hat{\xi}}^2(\bm{x}^*)$.
Since the inverse and determinant of the covariance matrix must be computed during GP training, the computational complexity is $O(N^3)$, where $N$ is the number of training samples. Therefore, the GP will encounter computational inefficiency on large datasets.

Note that the training data size for training $f_0$ and $f$ are $N_{\text{ED}}$ and $N_\text{ED} \cdot(n_W - 1)$, respectively. Generally speaking, for complex dynamical systems, $N_{\text{ED}}$ and $n_W$ are both about $O(10^1\sim10^2)$. Therefore, the GP can be directly used for emulating $f_0^{(1)},\cdots,f_0^{(m_{y})}$. However, when emulating $f^{(1)},\cdots,f^{(m_{y})}$, $N_\text{ED}\cdot (n_W - 1)$ (usually about $O(10^2\sim10^4)$) training samples typically constitute a large training set for GP, which will lead to efficiency issues. In this work, we utilize sparse GP (SGP) \cite{titsias2009variational} to address this problem.

Unlike conventional GP which utilizes the entire dataset for training, SGP adopts the concept of a sparse variational model. In this approach, a significantly smaller set of artificial training points, referred to as inducing points, is used to approximate the covariance structure of the full dataset, thereby greatly reducing computational cost while maintaining predictive accuracy. A trained SGP can also provide the predictive mean and variance of $\xi$ at any $\bm{x}^*$ as $\mu_{\hat{\xi}}(\bm{x}^*)$ and $\sigma_{\hat{\xi}}^2(\bm{x}^*)$. In SGP, only the covariance matrix associated with the inducing points needs to be inverted and its determinant evaluated. Thus, the computational complexity of SGP is $O(N_I^3)$, where $N_I$ is the number of inducing points. Since in most cases a few tens to a few hundreds of inducing points are sufficient, SGP is substantially more efficient than conventional GP when dealing with large training datasets.

Note that the time window width $T$ and the retained variance proportion $\varepsilon_{\lambda}$ for PCA are two key parameters in the F2NARX model, as they directly influence modeling accuracy. In this study, we employ cross-validation to determine suitable values of $T$ and $\varepsilon_{\lambda}$. Algorithm~\ref{alg:train-F2NARX} presents the pseudo-codes for the training approach of the F2NARX model with fixed $T$ and $\varepsilon_{\lambda}$.

\begin{algorithm}
\caption{Training approach of the F2NARX model}
\label{alg:train-F2NARX}
\KwIn{Training data $\left\{\bm{\theta}^{(i)},\bm{u}^{(i)},\bm{y}^{(i)}\right\}_{i=1}^{N_\text{ED}}$, time window width $T$, and $\varepsilon_{\lambda}$}
\KwOut{Surrogate model $\left\{\hat{\bm{f}}_0,\hat{\bm{f}}\right\}$ for emulating $y(t)$}

Obtain the number of time instants $n_T$ within each time window and the number of time windows $n_W$ for each training time histories according to $T$\;
Obtain matrices $\tilde{\bm{U}}$ and $\tilde{\bm{Y}}$ as shown in Eq.~(\ref{matrix for feature extraction})\;
Apply the feature extraction technique presented in Section~\ref{feature extraction} to obtain the feature mappings $\mathcal{K}_{u}$ and $\mathcal{K}_{y}$ from $\tilde{\bm{U}}$ and $\tilde{\bm{Y}}$, respectively, using the retained variance proportion $\varepsilon_{\lambda}$\;
$\bm{\xi}^{(i)}_{u,j} \gets \mathcal{K}_{u}(\tilde{\bm{u}}_j^{(i)}),\quad \bm{\xi}^{(i)}_{y,j} \gets \mathcal{K}_{y}(\tilde{\bm{y}}_j^{(i)}),\quad i=1,\cdots,N_\text{ED},\quad j=1,\cdots,n_W$\;
Generate the training data set $\{\bm{X}_0^\text{train},\bm{\varXi}_{y,0}^\text{train}\}$ as shown in Eq.~(\ref{training_data_f_0}) and use it to train the GP surrogate model $\hat{\bm{f}}_0$ for $\bm{f}_0$\;
Generate the training data set $\left\{\bm{X}^\text{train},\bm{\varXi}_{y}^\text{train}\right\}$ as shown in Eq.~(\ref{training_data_f}) and use it to train the SGP surrogate model $\hat{\bm{f}}$ for $\bm{f}$\;

\Return $\left\{\hat{\bm{f}}_0,\hat{\bm{f}}\right\}$\;

\end{algorithm}

\subsection{Probabilistic prediction of F2NARX model}

Because the surrogate models $\hat{\bm{f}}_0$ and $\hat{\bm{f}}$, whether GPs or SGPs, provide prediction uncertainties, this section explains how to leverage these uncertainties to generate probabilistic predictions of dynamical system responses.

To obtain the predicted mean of the response for a new sample of system parameters $\bm{\theta}^*$ and a new excitation $\bm{u}^*$, $\bm{u}^*$ is first segmented into $n_W$ local excitation vectors $\tilde{\bm{u}}^*_1,\cdots,\tilde{\bm{u}}^*_{n_W}$ by discarding the initial value. Then, $\hat{\bm{f}}_0(\cdot)$ is used to obtain the predictive mean of the low-dimensional features in the first time window as:
\begin{equation}
\bm{\mu}_{\hat{\bm{\xi}}^*_{y,1}}=\bm{\mu}_{\hat{\bm{f}}_0}\left(\left[\mathcal{K}_{u}(\tilde{\bm{u}}^*_1),u^*(t_0),y^*(t_0),\bm{\theta}^*\right]\right),
\end{equation}
where $\bm{\mu}_{\hat{\bm{f}}_0}(\cdot)=\left[\mu_{\hat{f}_0^{(1)}}(\cdot),\mu_{\hat{f}_0^{(2)}}(\cdot),\cdots,\mu_{\hat{f}_0^{(m_{y})}}(\cdot)\right]$ and $\mu_{\hat{f}_0^{(i)}}(\cdot)$ is the predictive mean function of the GP model $\hat{f}_0^{(i)}$; $u^*(t_0)$ and $y^*(t_0)$ represent the initial values of the excitation and response functions. Subsequently, the predictive mean of functional features in the subsequent time windows are predicted in an autoregressive manner as:
\begin{equation}
\bm{\mu}_{\hat{\bm{\xi}}^*_{y,j}}=\bm{\mu}_{\hat{\bm{f}}}\left(\left[\mathcal{K}_{u}(\tilde{\bm{u}}^*_j),\mathcal{K}_{u}(\tilde{\bm{u}}^*_{j-1}),\bm{\mu}_{\hat{\bm{\xi}}^*_{y,j-1}},\bm{\theta}^*\right]\right),j=2,\cdots,n_W,
\end{equation}
where $\bm{\mu}_{\hat{\bm{f}}}(\cdot)=\left[\mu_{\hat{f}^{(1)}}(\cdot),\mu_{\hat{f}^{(2)}}(\cdot),\cdots,\mu_{\hat{f}^{(m_{y})}}(\cdot)\right]$ and $\mu_{\hat{f}^{(i)}}(\cdot)$ is the predictive mean function of the SGP model $\hat{f}^{(i)}$. Then, the predicted mean of original dynamical response can be obtained as:
\begin{equation}
\bm{\mu}_{\hat{\bm{y}}^*}=\left[y^*(t_0),\mathcal{K}_{y}^{-1}\left(\bm{\mu}_{\hat{\bm{\xi}}^*_{{y},1}}\right),\cdots,\mathcal{K}_{y}^{-1}\left(\bm{\mu}_{\hat{\bm{\xi}}^*_{y,n_W}}\right)\right].
\end{equation}

\begin{figure}[t]
\centering
\includegraphics[scale=0.45]{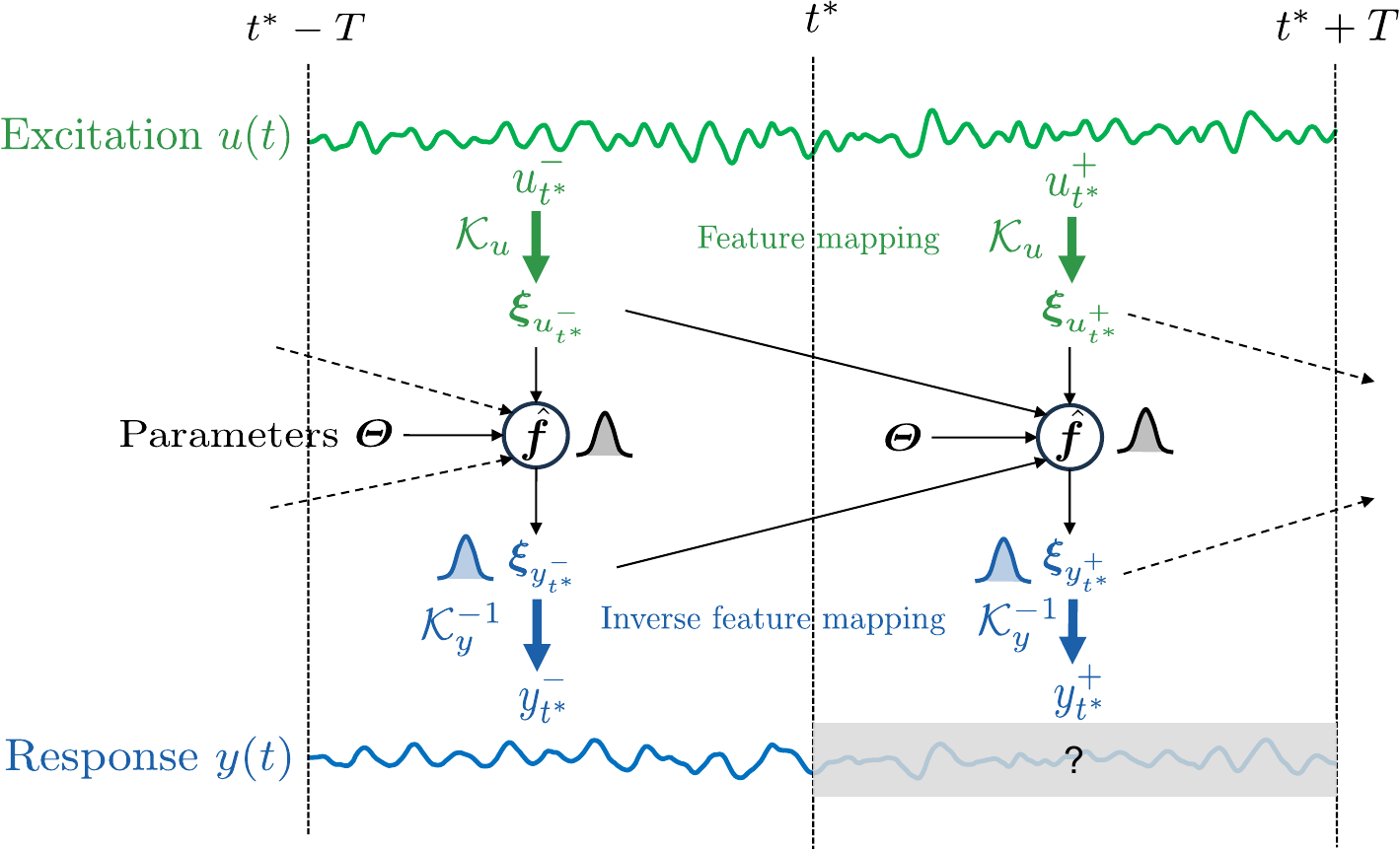}
\caption{Illustration of probabilistic prediction approach of the F2NARX model.}
\label{Prob_pred_F2NARX}
\end{figure}

Estimating the predictive variance of the response is more challenging, as F2NARX employs an autoregressive structure in which each time window's prediction will be as input for the next time window. Consequently, both the surrogate modeling uncertainty and the uncertainty propagated from earlier predicted response must be taken into account. Fig.~\ref{Prob_pred_F2NARX} gives an illustration for this.

For the first time window, response uncertainty arises solely from the surrogate modeling uncertainty; hence, the prediction covariance matrix of the functional features is:
\begin{equation}
    \bm{\varSigma}_{\hat{\bm{\xi}}^*_{y,1}}=\bm{\varSigma}_{\hat{\bm{f}}_0}(\left[\mathcal{K}_{u}(\tilde{\bm{u}}^*_1),u^*(t_0),y^*(t_0),\bm{\theta}^*\right]),
\end{equation}
where $\bm{\varSigma}_{\hat{\bm{f}}_0}(\cdot)=\text{diag}\left(\left[\sigma^2_{\hat{f}_0^{(1)}}(\cdot),\sigma^2_{\hat{f}_0^{(2)}}(\cdot),\cdots,\sigma^2_{\hat{f}_0^{(m_{y})}}(\cdot)\right]\right)$ and $\sigma^2_{\hat{f}_0^{(i)}}(\cdot)$ is the predictive variance function of the GP model $\hat{f}_0^{(i)}$. For subsequent time windows, response uncertainty stems from both surrogate modeling uncertainty and the propagated uncertainty of the preceding window’s predictions. Because the low-dimensional features are uncorrelated, we first show how to propagate both sources of uncertainty for a single feature. Without loss of generality, let $\xi(\bm{x})$ denote the SGP model for a single feature and let $\bm{x}^*$ be the input vector such that $\bm{x}^* \sim \mathcal{N}(\bm{\mu}_{\bm{x}^*}, \bm{\varSigma}_{\bm{x}^*})$, where $\bm{\mu}_{\bm{x}^*}$ and $\bm{\varSigma}_{\bm{x}^*}$ are the input mean vector and covariance matrix, respectively. Then, the problem becomes how to estimate $\mathrm{Var}(\xi(\bm{x}^*))$ while accounting for uncertainty in both $\bm{x}^*$ and $\xi(\cdot)$.

A straightforward approach is to use MCS. 
However, this requires a large number of samples and is highly inefficient for autoregressive models. To estimate the prediction variance efficiently for autoregressive models, Girard et al. \cite{girard2002multiple} developed an analytical expression of $\mathrm{Var}(\xi(\bm{x}^*))$ as:
\begin{equation}
\label{eq_pred_variance}
    \mathrm{Var}(\xi(\bm{x}^*))=\mathbb{E}_{\bm{x}^*}[\sigma^2_{\xi}(\bm{x}^*)]+\mathrm{Var}_{\bm{x}^*}(\mu_{\xi}(\bm{x}^*)),
\end{equation}
where $\mu_{\xi}(\bm{x}^*)$ and $\sigma^2_{\xi}(\bm{x}^*)$ are the predicted mean and variance of SGP model. Based on this, they proposed a Taylor expansion-based approach to estimate $\mathrm{Var}(\xi(\bm{x}^*))$ efficiently, where they approximate $\mu_{\xi}(\bm{x}^*)$ and $\sigma^2_{\xi}(\bm{x}^*)$ by their first order and second order Taylor expansion around $\mu_{\bm{x}^*}$ respectively. 

Here, we propose an alternative, more efficient approach to estimate $\mathrm{Var}(\xi(\bm{x}^*))$. Eq.~(\ref{eq_pred_variance}) indicates that we only need to compute the mean of the variance function $\sigma^2_{\xi}(\cdot)$ and the variance of the mean function $\mu_{\xi}(\cdot)$ under the distribution of $\bm{x}^*$. We propose to use the unscented transform to estimate them. The unscented transform (UT) \cite{julier1997new} is an efficient method to approximate the statistics (mean and variance) of a nonlinear transformation $h(\cdot)$ of a random variable. It picks a small, carefully chosen set of sample points, called sigma points, around the mean value of $\bm{x}^*$ so that they capture the first two moments (mean and covariance) of $\bm{x}^*$. Then it estimates the mean and variance of the nonlinear transformation of $\bm{x}^*$ through weighted sample mean and variance. For a $p$-dimensional $\bm{x}^*$, the UT typically uses $2p+1$ sigma points associated with weights as:
\begin{equation}
\begin{gathered}
    \bm{s}_0=\bm{\mu}_{\bm{x}^*},\quad \bm{s}_i=\bm{\mu}_{\bm{x}^*}+\left(\sqrt{(p+\kappa)\bm{\varSigma}_{\bm{x}^*}}\right)_i,\quad
    \bm{s}_{i+p}=\bm{\mu}_{\bm{x}^*}-\left(\sqrt{(p+\kappa)\bm{\varSigma}_{\bm{x}^*}}\right)_i,\\
    \alpha_0=\frac{\kappa}{p+\kappa},\quad 
    \alpha_i=\frac{1}{2(p+\kappa)},\quad
    \alpha_{i+p}=\frac{1}{2(p+\kappa)} \quad i=1,\cdots,p,
\end{gathered}
\end{equation}
where $\left(\sqrt{(p+\kappa)\bm{\varSigma}_{\bm{x}^*}}\right)_i$ is the $i$-th row of the matrix square root of $(p+\kappa)\bm{\varSigma}_{\bm{x}^*}$ and $\kappa$ is a tuning parameter.
As a heuristic for Gaussian variables, $\kappa$ is selected so that $p+\kappa=3$ \cite{julier1997new}.
Note that $\kappa$ can be either positive or negative according to $p$.
For high-dimensional effective inputs, scaled unscented transform \cite{julier2002scaled} could be used.
Then, the mean and variance of $h(\bm{x}^*)$ are given by \cite{julier1997new}:
\begin{equation}
\label{eq_mean_UT}
    \mathbb{E}_{\bm{x}^*}[h(\bm{x}^*)] \approx \sum_{i=0}^{2p} \alpha_i h(\bm{s}_i),
\end{equation}
\begin{equation}
\label{eq_var_UT}
\begin{split}
    \mathrm{Var}_{\bm{x}^*}(h(\bm{x}^*)) 
    &\approx \sum_{i=0}^{2p} \alpha_i (h(\bm{s}_i)-\mathbb{E}_{\bm{x}^*}[h(\bm{x}^*)])^2\\
    &=\sum_{i=0}^{2p} \alpha_i h^2(\bm{s}_i)-\left[\sum_{i=0}^{2p} \alpha_i h(\bm{s}_i)\right]^2.
\end{split}
\end{equation}

By substituting $h(\cdot)$ with the predicted mean function $\mu_{\xi}(\cdot)$ and predicted variance function $\sigma^2_{\xi}(\cdot)$ in Eq.~\eqref{eq_var_UT} and Eq.~\eqref{eq_mean_UT} respectively, we can obtain the variance of $\mu_{\xi}(\bm{x}^*)$ and mean of $\sigma^2_{\xi}(\bm{x}^*)$ as:
\begin{equation}
\label{eq_UT_var_of_mean}
    \mathrm{Var}_{\bm{x}^*}(\mu_{\xi}(\bm{x}^*)) 
    \approx
    \sum_{i=0}^{2p} \alpha_i \mu_{\xi}^2(\bm{s}_i)-\left[\sum_{i=0}^{2p} \alpha_i \mu_{\xi}(\bm{s}_i)\right]^2,
\end{equation}
\begin{equation}
\label{eq_UT_mean_of_var}
    \mathbb{E}_{\bm{x}^*}[\sigma^2_{\xi}(\bm{x}^*)] \approx \sum_{i=0}^{2p} \alpha_i \sigma^2_{\xi}(\bm{s}_i).
\end{equation}
Then, by substituting Eq.~(\ref{eq_UT_var_of_mean}) and Eq.~(\ref{eq_UT_mean_of_var}) into Eq.~(\ref{eq_pred_variance}) we obtain the approximate predicted variance:
\begin{equation}
\label{pred_var_UT}
    \mathrm{Var}(\xi(\bm{x}^*)) 
    \approx \sum_{i=0}^{2p} \alpha_i \left[\sigma^2_{\xi}(\bm{s}_i)+\mu_{\xi}^2(\bm{s}_i)\right]-\left[\sum_{i=0}^{2p} \alpha_i \mu_{\xi}(\bm{s}_i)\right]^2.
\end{equation}

For the subsequent windows, the predictive covariance matrices of the low-dimensional features are:
\begin{equation}
    \bm{\varSigma}_{\hat{\bm{\xi}}^*_{y,j}}=\text{diag}\left(\left[\mathrm{Var}\left(\hat{f}^{(1)}(\bm{x}^*_j)\right),\mathrm{Var}\left(\hat{f}^{(2)}(\bm{x}^*_j)\right),\cdots,\mathrm{Var}\left(\hat{f}^{(m_{y})}(\bm{x}^*_j)\right)\right]\right),j=2,\cdots,n_W,
\end{equation}
where each $\mathrm{Var}\left(\hat{f}^{(i)}(\bm{x}^*_j)\right)$ is obtained from Eq.~(\ref{pred_var_UT}) by substituting $\xi(\cdot)$ with $\hat{f}^{(i)}(\cdot)$ for $i=1,2,\cdots,m_{y}$. Here, $\bm{x}^*_j \sim \mathcal{N}(\bm{\mu}_{\bm{x}^*_j}, \bm{\varSigma}_{\bm{x}^*_j})$ with
\begin{equation}
\label{input_mean_cov}
\begin{gathered}
    \bm{\mu}_{\bm{x}^*_j}=\left[\mathcal{K}_{u}(\tilde{\bm{u}}^*_j),\mathcal{K}_{u}(\tilde{\bm{u}}^*_{j-1}),\bm{\mu}_{\hat{\bm{\xi}}^*_{y,j-1}},\bm{\theta}^*\right],\\
    \bm{\varSigma}_{\bm{x}^*_j}=\text{diag}\left(\left[\bm{0}_{1\times m_u},\bm{0}_{1\times m_u},\text{diag}\left(\bm{\varSigma}_{\hat{\bm{\xi}}^*_{y,j-1}}\right),\bm{0}_{1\times n_s}\right]\right).
\end{gathered}
\end{equation}
Note that the operator $\text{diag}(\cdot)$ behaves differently depending on its input: if the input is a vector, the output is a square diagonal matrix with the vector elements on the main diagonal; if the input is a square matrix, the output is a vector containing the elements from the main diagonal of the input matrix. According to Eq.~(\ref{inverse-feature-mapping-expression}), the prediction variance of original dynamical response can be given as:
\begin{equation}
\bm{\sigma}^2_{\hat{\bm{y}}^*}=\left[0,\text{diag}\left(\bm{V}_{y}\bm{\varSigma}_{\hat{\bm{\xi}}^*_{y,1}}\bm{V}_{y}^\top\right),\cdots,\text{diag}\left(\bm{V}_{y}\bm{\varSigma}_{\hat{\bm{\xi}}^*_{y,n_W}}\bm{V}_{y}^\top\right)\right].
\end{equation}
Algorithm~\ref{alg:predict-F-NARX-F} presents the pseudo-codes for the probabilistic predicting approach of the F2NARX model.

\begin{algorithm}
\caption{Predicting approach of the F2NARX model}
\label{alg:predict-F-NARX-F}
\KwIn{Surrogate model $\left\{\hat{\bm{f}}_0,\hat{\bm{f}}\right\}$ for emulating $y(t)$, a new sample of parameters $\bm{\theta}^*$, and a new excitation $\bm{u}^*$}
\KwOut{Predictive mean $\bm{\mu}_{\hat{\bm{y}}^*}$ and variance $\bm{\sigma}^2_{\hat{\bm{y}}^*}$ of dynamical system response}

Segment $\bm{u}^*$ into $n_W$ local excitation vectors $\tilde{\bm{u}}^*_1,\cdots,\tilde{\bm{u}}^*_{n_W}$ by discarding the initial value\;

$\bm{\mu}_{\hat{\bm{\xi}}^*_{y,1}} \gets \bm{\mu}_{\hat{\bm{f}}_0}\left(\left[\mathcal{K}_{u}(\tilde{\bm{u}}^*_1),u^*(t_0),y^*(t_0),\bm{\theta}^*\right]\right),\quad \bm{\varSigma}_{\hat{\bm{\xi}}^*_{y,1}} \gets \bm{\varSigma}_{\hat{\bm{f}}_0}(\left[\mathcal{K}_{u}(\tilde{\bm{u}}^*_1),u^*(t_0),y^*(t_0),\bm{\theta}^*\right])$\;

\For{$j = 2$ \KwTo $n_W$}{
    $\bm{\mu}_{\hat{\bm{\xi}}^*_{y,j}} \gets \bm{\mu}_{\hat{\bm{f}}}\left(\left[\mathcal{K}_{u}(\tilde{\bm{u}}^*_j),\mathcal{K}_{u}(\tilde{\bm{u}}^*_{j-1}),\bm{\mu}_{\hat{\bm{\xi}}^*_{y,j-1}},\bm{\theta}^*\right]\right)$\;

    $\bm{\varSigma}_{\hat{\bm{\xi}}^*_{y,j}} \gets \text{diag}\left(\left[\mathrm{Var}\left(\hat{f}^{(1)}(\bm{x}^*_j)\right),\mathrm{Var}\left(\hat{f}^{(2)}(\bm{x}^*_j)\right),\cdots,\mathrm{Var}\left(\hat{f}^{(m_{y})}(\bm{x}^*_j)\right)\right]\right)$ with $\bm{x}^*_j \sim \mathcal{N}(\bm{\mu}_{\bm{x}^*_j}, \bm{\varSigma}_{\bm{x}^*_j})$, where $\bm{\mu}_{\bm{x}^*_j}$ and $\bm{\varSigma}_{\bm{x}^*_j}$ are obtained through Eq.~(\ref{input_mean_cov})\;
}

$\bm{\mu}_{\hat{\bm{y}}^*} \gets \left[y^*(t_0),\mathcal{K}_{y}^{-1}\left(\bm{\mu}_{\hat{\bm{\xi}}^*_{{y},1}}\right),\cdots,\mathcal{K}_{y}^{-1}\left(\bm{\mu}_{\hat{\bm{\xi}}^*_{y,n_W}}\right)\right]$\;

$\bm{\sigma}^2_{\hat{\bm{y}}^*} \gets \left[0,\text{diag}\left(\bm{V}_{y}\bm{\varSigma}_{\hat{\bm{\xi}}^*_{y,1}}\bm{V}_{y}^\top\right),\cdots,\text{diag}\left(\bm{V}_{y}\bm{\varSigma}_{\hat{\bm{\xi}}^*_{y,n_W}}\bm{V}_{y}^\top\right)\right]$\;

\Return $\bm{\mu}_{\hat{\bm{y}}^*}$ and $\bm{\sigma}^2_{\hat{\bm{y}}^*}$\;

\end{algorithm}

\section{Active learning with probabilistic F2NARX predictions for first-passage probability estimation}
\label{sec_AL}

To demonstrate the usefulness of the probabilistic prediction capability of F2NARX, we employ the proposed model to guide active learning for first-passage probability estimation of stochastic structural dynamical systems. 
This section briefly introduces the first-passage probability and the adopted active learning strategy. 
It should be emphasized that the purpose of this section is not to develop a general method for rare-event reliability analysis, but rather to illustrate the effectiveness of probabilistic F2NARX and its potential for reliability-oriented tasks.

\subsection{First-passage probability of structural dynamical systems}
With a structural dynamical system described by Eq.~(\ref{governing_eq_MDOF}), the dynamical response of interest $y(\bm{\varTheta},\bm{\varPhi},t)$ depends on both the uncertain parameter vector $\bm{\varTheta}$ and the random vector $\bm{\varPhi}$ controlling the randomness in excitation. 
We consider a double-sided boundary condition to define failure occurrence. The first-passage probability is defined as the probability that $|y(\bm{\varTheta},\bm{\varPhi},t)|$ exceeds a prescribed safety threshold $y_{\text{th}}$ for the first time within the interval $[t_0,t_e]$, which is expressed as:
\begin{equation}
    \label{definition_of_Pf}
    P_f=\text{Pr}\{|y(\bm{\varTheta},\bm{\varPhi},t)|\geq y_{\text{th}},\exists t\in[t_0,t_e]\},
\end{equation}
where $\text{Pr}\{\cdot\}$ is the probability operator, $\exists$ means `there exists'. It can be further expressed as:
\begin{equation}
    \label{definition_of_FPP_integral}
    \begin{split}
    P_f &= \text{Pr}\left\{\max_{t\in[t_0,t_e]}|y(\bm{\varTheta},\bm{\varPhi},t)|\geq y_{\text{th}}\right\}  \\
    &=\int_{\varOmega_{\bm{\varTheta}}}\int_{\varOmega_{\bm{\varPhi}}}I\left(\max_{t\in[t_0,t_e]}|y(\bm{\theta},\bm{\phi},t)|\geq y_{\text{th}}\right)f_{\bm{\varTheta}}(\bm{\theta})f_{\bm{\varPhi}}(\bm{\phi})\text{d}\bm{\theta}\text{d}\bm{\phi},
    \end{split}
\end{equation}
where $I(\cdot)$ is the indicator function, $f_{\bm{\varTheta}}(\bm{\theta})$ is the joint probability density function of $\bm{\varTheta}$, and $f_{\bm{\varPhi}}(\bm{\phi})$ is the joint probability density function of $\bm{\varPhi}$.

In general, a closed-form solution of Eq.~(\ref{definition_of_FPP_integral}) is not available. Therefore, MCS is typically employed to estimate $P_f$. Given an MCS sample set $\left\{\bm{\theta}^{(i)}\right\}_{i=1}^{N_{\text{MCS}}}$ and $\left\{\bm{\phi}^{(i)}\right\}_{i=1}^{N_{\text{MCS}}}$ of $\bm{\varTheta}$ and $\bm{\varPhi}$, drawn from $f_{\bm{\varTheta}}(\bm{\theta})$ and $f_{\bm{\varPhi}}(\bm{\phi})$, respectively, the dynamical system is evaluated at each $\left\{\bm{\theta}^{(i)},\bm{\phi}^{(i)}\right\}$ to obtain the QoI $y^{(i)}(t)=y(\bm{\theta}^{(i)},\bm{\phi}^{(i)},t)$.
The MCS estimator of the first-passage probability is then given by:
\begin{equation}
    \label{MCS_estimation_Pf}
    \hat{P}_f=\frac{1}{N_{\text{MCS}}}\sum_{i=1}^{N_{\text{MCS}}}I\left(\max_{j\in\{1,\cdots,N_t\}}|y^{(i)}(t_j)|\geq y_{\text{th}}\right),
\end{equation}
where $\{t_1=t_0,t_2,\cdots,t_{N_t}=t_e\}$ is an evenly spaced grid over $[t_0,t_e]$.

Evaluating Eq.~(\ref{MCS_estimation_Pf}) requires a large number of dynamical system simulations, which becomes computationally prohibitive when the system is expensive to simulate. To alleviate this, we construct a cheap-to-evaluate surrogate model $\hat{y}(\bm{\varTheta},\bm{\varPhi},t)$ using the proposed F2NARX method. For the first-passage probability estimation considered here, the surrogate model primarily needs to be accurate near the limit-state surface, i.e., the boundary between the safe and failure domains, rather than over the entire input space. Therefore, training the F2NARX model using samples close to the limit-state surface can further reduce the required training dataset size. Active learning \cite{moustapha2022active} provides an effective way to adaptively select such informative samples.

\subsection{Active learning strategy}
\begin{figure}[t]
\centering
\includegraphics[scale=0.6]{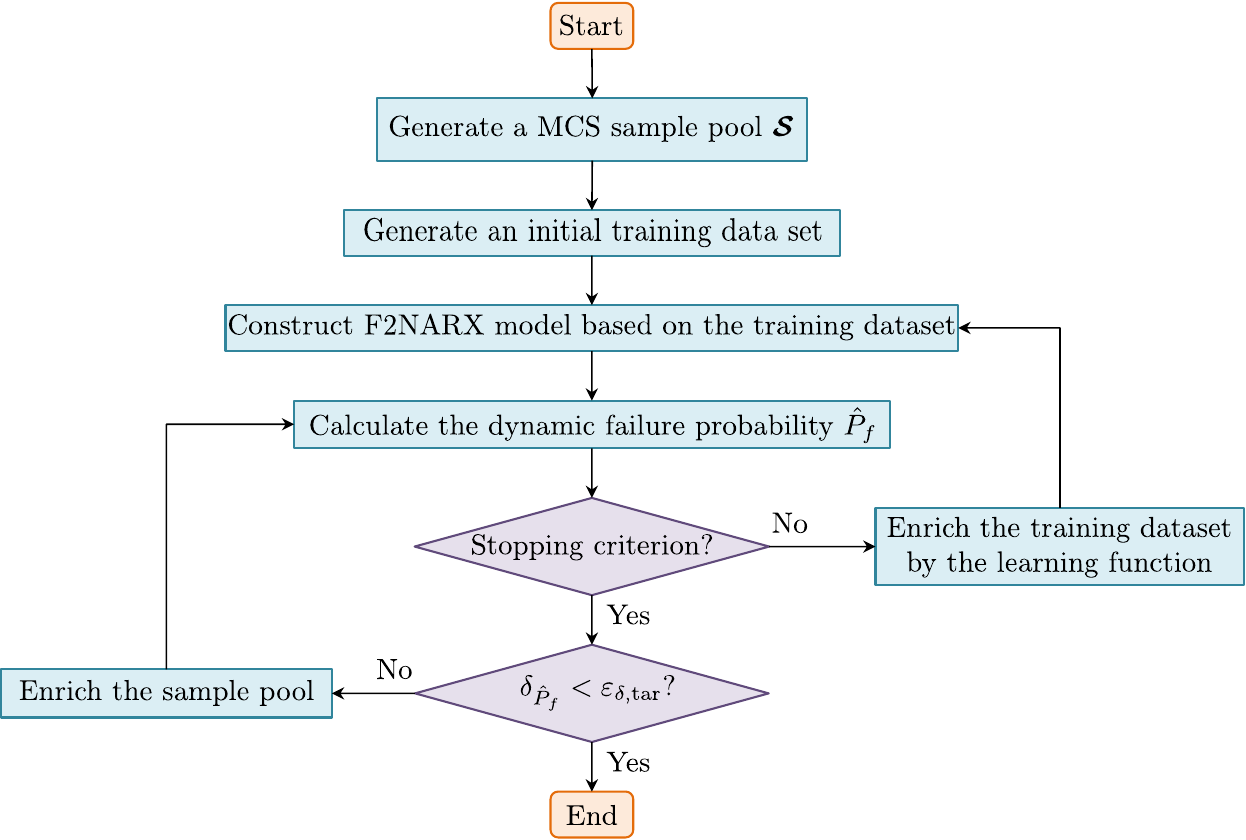}
\caption{Flowchart of the active learning procedure for first-passage probability estimation.}
\label{flowchart_AL}
\end{figure}

Active learning begins with an initial training set and sequentially enriches it by selecting the most informative samples according to specific criteria. Two key components of active learning are the learning function and the stopping criterion \cite{moustapha2022active}. In this work, we adapt the learning function from the well-known time-dependent reliability analysis method SILK \cite{hu2016single} and modify it to guide the selection of the most informative next time history. For stopping criterion, a sampling resource-based criterion is employed to decide when the active learning process should stop. The overall flowchart of the active learning procedure is shown in Fig.~\ref{flowchart_AL}, and the detailed implementation steps are summarized below.

\textbf{Step 1: Generate a sample pool}

Generate an MCS sample pool $\bm{\mathcal{S}}=\left\{\left(\bm{\theta}^{(i)},\bm{\phi}^{(i)}\right)\right\}_{i=1}^{N_{\text{MCS}}}$ according to $f_{\bm{\varTheta}}(\bm{\theta})$ and $f_{\bm{\varPhi}}(\bm{\phi})$.

\textbf{Step 2: Generate an initial training dataset}

Generate an initial dataset $\left\{\left(\bm{\theta}^{(i)},\bm{\phi}^{(i)}\right)\right\}_{i=1}^{N}$ using MCS. Construct the excitation dataset $\left\{\bm{u}^{(i)}\right\}_{i=1}^{N}$ from it, and obtain the dynamical response dataset $\left\{\bm{y}^{(i)}\right\}_{i=1}^{N}$. Let $N_{\text{new}}$ denote the number of training samples added during the active learning process, and initialize $N_{\text{new}}=0$.

\textbf{Step 3: Construct the surrogate model}

Construct the F2NARX model to emulate $y(\bm{\varTheta},\bm{\varPhi},t)$ based on the training dataset $\left\{\bm{\theta}^{(i)},\bm{u}^{(i)},\bm{y}^{(i)}\right\}_{i=1}^{N}$ according to Algorithm \ref{alg:train-F2NARX}.

\textbf{Step 4: Estimate the first-passage failure probability}

Use Algorithm \ref{alg:predict-F-NARX-F} to obtain probabilistic predictions with the input $\bm{\mathcal{S}}$ and use the prediction mean to calculate the first-passage failure probability $\hat{P}_f$.

\textbf{Step 5: Check the stopping criterion}

A sampling resource-based stopping criterion is used here. Let $N^{\text{target}}_{\text{new}}$ denote the given sampling resource. If $N_{\text{new}}=N^{\text{target}}_{\text{new}}$ is satisfied, proceed to \textbf{Step 7}; otherwise, go to \textbf{Step 6}.

\textbf{Step 6: Enrich the training dataset by the learning function}

Here, we utilize the double-sided $U_{\text{min}}$ function to select the best next point. For a sample $\left(\bm{\theta}^*,\bm{\phi}^*\right)\in\mathcal{S}$, the double-sided $U_{\text{min}}$ is defined as:
\begin{equation}
    U_{\text{min}}(\bm{\theta}^*,\bm{\phi}^*)=\min\left(U_{\text{min}}^{\text{upper}}(\bm{\theta}^*,\bm{\phi}^*),U_{\text{min}}^{\text{lower}}(\bm{\theta}^*,\bm{\phi}^*)\right),
\end{equation}
where $U_{\text{min}}^{\text{upper}}(\bm{\theta}^*,\bm{\phi}^*)$ and $U_{\text{min}}^{\text{lower}}(\bm{\theta}^*,\bm{\phi}^*)$ are single-sided $U_{\text{min}}$ functions calculated as \cite{hu2016single}:
\begin{equation}
    U_{\text{min}}^{\text{upper}}(\bm{\theta}^*,\bm{\phi}^*)=
    \begin{cases}
    u_e,\ \exists j\in\{1,\cdots,N_t\},\mu_{\hat{y}}(\bm{\theta}^*,\bm{\phi}^*,t_j)\geq y_{\text{th}}\ \text{and}\ U^{\text{upper}}(\bm{\theta}^*,\bm{\phi}^*,t_j)\geq2, \\
    \min\limits_{j=1,\dots,N_t} U^{\text{upper}}(\bm{\theta}^*,\bm{\phi}^*,t_j),\ \text{otherwise},
    \end{cases}
\end{equation}
\begin{equation}
    U_{\text{min}}^{\text{lower}}(\bm{\theta}^*,\bm{\phi}^*)=
    \begin{cases}
    u_e,\ \exists j\in\{1,\cdots,N_t\},\mu_{\hat{y}}(\bm{\theta}^*,\bm{\phi}^*,t_j)\leq-y_{\text{th}}\ \text{and}\ U^{\text{lower}}(\bm{\theta}^*,\bm{\phi}^*,t_j)\geq2, \\
    \min\limits_{j=1,\dots,N_t} U^{\text{lower}}(\bm{\theta}^*,\bm{\phi}^*,t_j),\ \text{otherwise},
    \end{cases}
\end{equation}
where $u_e$ is a value larger than two \cite{hu2016single} and
\begin{equation}
    U^{\text{upper}}(\bm{\theta}^*,\bm{\phi}^*,t_j)=\frac{|\mu_{\hat{y}}(\bm{\theta}^*,\bm{\phi}^*,t_j)-y_{\text{th}}|}{\sigma_{\hat{y}}(\bm{\theta}^*,\bm{\phi}^*,t_j)},
\end{equation}
\begin{equation}
    U^{\text{lower}}(\bm{\theta}^*,\bm{\phi}^*,t_j)=\frac{|\mu_{\hat{y}}(\bm{\theta}^*,\bm{\phi}^*,t_j)+y_{\text{th}}|}{\sigma_{\hat{y}}(\bm{\theta}^*,\bm{\phi}^*,t_j)},
\end{equation}
in which $\mu_{\hat{y}}(\bm{\theta}^*,\bm{\phi}^*,t_j)$ and $\sigma_{\hat{y}}(\bm{\theta}^*,\bm{\phi}^*,t_j)$ are the predicted mean and standard deviation at $(\bm{\theta}^*,\bm{\phi}^*,t_j)$ provided by F2NARX. The next sample point $(\bm{\theta}^{\text{new}},\bm{\phi}^{\text{new}})$ is then selected by minimizing the double-sided $U_{\text{min}}$ function as:
\begin{equation}
    \left(\bm{\theta}^{\text{new}},\bm{\phi}^{\text{new}}\right)=\mathop{\mathrm{argmin}}\limits_{\left(\bm{\theta},\bm{\phi}\right)\in\bm{\mathcal{S}}}\ U_{\text{min}}(\bm{\theta},\bm{\phi}).
\end{equation}
Next, enrich the current training dataset with the new sample $\left\{\bm{\theta}^{\text{new}},\bm{u}^{\text{new}},\bm{y}^{\text{new}}\right\}$ and update $N_{\text{new}} = N_{\text{new}}+1$. Then proceed to \textbf{Step 3}.

\textbf{Step 7: Check the CoV of the first-passage probability estimate}

Calculate the CoV of $\hat{P}_f$ as follows:
\begin{equation}
    \delta_{\hat{P}_f}=\sqrt{\frac{1-\hat{P}_f}{(N_{\text{MCS}}-1)\hat{P}_f}}.
\end{equation}
If $\delta_{\hat{P}_f}<\varepsilon_{\delta,\text{tar}}$ is satisfied, the active learning approach \textbf{ends}; otherwise, go to \textbf{Step 8}.
Note that this CoV check is not a convergence criterion for the surrogate quality but is used to control the sampling error associated with the MCS sample pool.

\textbf{Step 8: Enrich the sample pool}

Enrich the current sample pool by adding new samples, i.e., $\bm{\mathcal{S}} = \bm{\mathcal{S}} \cup \bm{\mathcal{S}}^+$, where $\bm{\mathcal{S}}^+$ is generated in the same manner as in \textbf{Step 1}. Then, proceed to \textbf{Step 4}. This step aims to enlarge the MCS sample pool if it is insufficient to estimate the failure probability with the desired CoV.

\section{Examples and discussions}
\label{sec_4}
Two case studies, a single-degree-of-freedom Bouc-Wen oscillator and a realistic nonlinear three-story steel frame, are presented in this section to evaluate the performance of the proposed F2NARX method. 
All computations are performed on a MacBook Pro equipped with an Apple M4 Pro chip and 24 GB of RAM.
All GP and SGP models are trained using the GPyTorch Python package \cite{gardner2018gpytorch}.
SGP models use 500 inducing points, whose locations are treated as model hyperparameters. 
The hyperparameters of both GP and SGP models are optimized using the Adam optimizer in PyTorch with a learning rate of 0.05 and 200 optimization iterations.
In this work, the CPU version of PyTorch was used for GPyTorch-based GP and SGP training and prediction. 
In addition, all GP and SGP models in F2NARX were trained and used for prediction sequentially rather than in parallel.
For each example, results are reported from the following four aspects: 

\textbf{(1) Surrogate modeling error of the F2NARX method under various parameter settings}, including time window width and retained variance proportion for PCA. The mean prediction error over a testing dataset is given by:
\begin{equation}
\label{pred_error}
    \bar{\epsilon}=\frac{1}{N_{\text{test}}}\sum_{i=1}^{N_{\text{test}}}\epsilon\left(\bm{y}^{(i)},\hat{\bm{y}}^{(i)}\right),
\end{equation}
where $N_{\text{test}}$ is the size of testing set and $\epsilon\left(\bm{y}^{(i)},\hat{\bm{y}}^{(i)}\right)$ is the normalized mean squared error between the true response $\bm{y}^{(i)}$ and predicted response $\hat{\bm{y}}^{(i)}$:
\begin{equation}
\label{NMSE}
    \epsilon\left(\bm{y}^{(i)},\hat{\bm{y}}^{(i)}\right)=\frac{1}{N_t}\frac{\left\|\bm{y}^{(i)}-\hat{\bm{y}}^{(i)}\right\|_2^2}{\text{Var}\left(\bm{y}^{(i)}\right)},
\end{equation}
where $N_t$ is the length of the response vector and $\|\cdot\|_2$ is the Euclidean norm.

\textbf{(2) Comparison of F2NARX with an SGP-based variant of the state-of-the-art $\mathcal{F}$-NARX model} in terms of modeling error, training time, and prediction time under different training data sizes.
The SGP-based variant of $\mathcal{F}$-NARX, abbreviated as SGP-F-NARX, adopts the same formulation and feature extraction method as the original $\mathcal{F}$-NARX model \cite{schar2025surrogate}.
The difference between SGP-F-NARX and $\mathcal{F}$-NARX lies in the surrogate modeling of $\mathcal{F}(\cdot)$ in Eq.~(\ref{F-NARX model}): instead of the polynomial regression adopted in \cite{schar2025surrogate}, SGP-F-NARX employs an SGP model to compare with F2NARX.
In this way, the comparison mainly reflects the difference between the one-step-ahead autoregressive structure of $\mathcal{F}$-NARX and the one-time-window-ahead autoregressive structure of F2NARX, rather than the difference between different regression models.

\textbf{(3) Modeling error and prediction time for probabilistic predictions} using the F2NARX method. Two probabilistic prediction methods are considered: the Taylor expansion-based method proposed in \cite{girard2002multiple} and the unscented transform-based method proposed in this study. The mean prediction errors (computed in the same way as Eq.~(\ref{pred_error})) of the prediction standard deviation function obtained by the Taylor expansion-based and unscented transform-based methods, relative to the MCS results, are reported, along with their corresponding computational times.

\textbf{(4) Active learning results} based on the probabilistic predictions for estimating the first-passage failure probability of dynamical systems, along with a comparison against results obtained without using active learning. 
The target CoV of the failure probability estimate is set to $0.05$ in the examples.
Since the purpose of this comparison is to demonstrate the effectiveness of probabilistic F2NARX and its potential for reliability-oriented tasks, rather than to evaluate F2NARX as a general rare-event reliability analysis method, the examples consider moderate failure probabilities, specifically $0.0405$ and $0.0411$.

\subsection{Example 1: The single-degree-of-freedom Bouc-Wen oscillator}
\begin{figure}[t]
\centering
\includegraphics[scale=0.4]{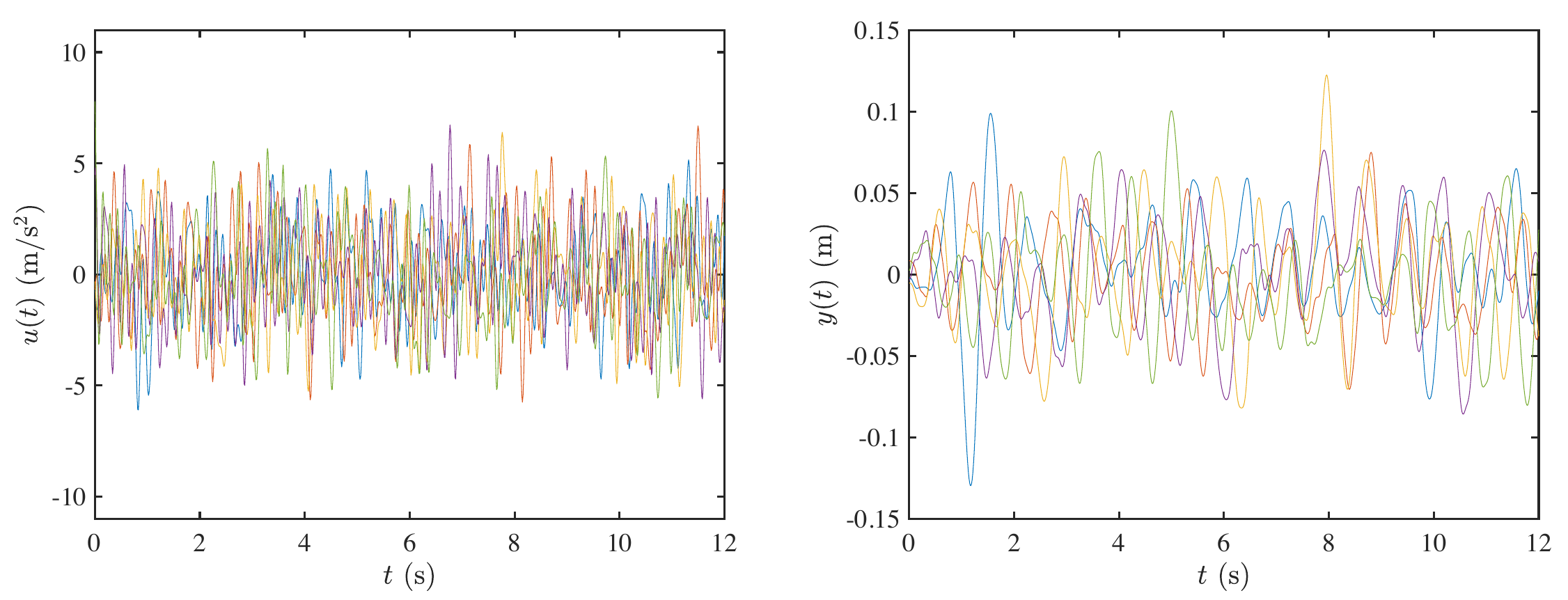}
\caption{Five different realizations of excitations $u(t)$ and corresponding responses $y(t)$ for the Bouc-Wen oscillator.}
\label{Ex1_u(t)_y(t)}
\end{figure}

The first example investigates a Bouc-Wen oscillator \cite{mai2016surrogate,papaioannou2019improved,cheng2025state} governed by the following ordinary differential equations:
\begin{equation}
    \label{eq_Bouc_Wen}
    \begin{cases}
m\ddot{y}(t)+c\dot{y}(t)+k\left[\alpha y(t)+(1-\alpha)x_yz(t)\right]=mu(t), \\
\dot{z}(t) = \frac{1}{x_y}\left[A\dot{y}(t)-\beta \left|\dot{y}(t)\right|\left|z(t)\right|^{n-1}z(t)-\gamma \dot{y}(t)\left|z(t)\right|^{n}\right],
\end{cases}
\end{equation}
where $y(t)$ is the oscillator displacement with initial conditions $\dot{y}(0)=0$ and $y(0)=y_0$, and $z(t)$ is the hysteretic displacement with initial condition $z(0)=0$.  $m$, $k$, and $y_0$ are treated as independent uniform random variables such that $m\sim U(5\times10^4,7\times10^4)\ \text{kg}$, $k\sim U(4\times10^6,6\times10^6)\ \text{N}/\text{m}$, and $y_0\sim U(-1\times10^{-2},1\times10^{-2})\ \text{m}$. The excitation $u(t)$ is modeled through spectral representation as:
\begin{equation}
    u(t)=S_u\sum_{i=1}^{500}[\phi_{i}\cos(\omega_it)+\phi_{i+500}\sin(\omega_it)],
\end{equation}
where $\Delta\omega=\omega_u/500$, $\omega_u$ is the upper cut-off frequency and is set to $15\pi~\mathrm{rad/s}$, $\omega_i=i\Delta\omega$, and $\phi_{i},i=1,2,\cdots,1000$ are independent standard Gaussian random variables. In addition, we set $c=0.1m\sqrt{k/m}$, $\alpha=\beta=\gamma=0.5$, $A=1$, $n=3$, $x_y=0.04\  \text{m}$, and $S_u=9.70813\times10^{-2}~\text{m}/\text{s}^2$. The quantity of interest is $y(t)$ over the time interval [0,12] s and the time interval is discretized into 3001 equally spaced time instants with time increment $\delta t=0.004\ \text{s}$. Fig.~\ref{Ex1_u(t)_y(t)} shows five different realizations for this problem.

\subsubsection{Performance of F2NARX method under different parameter settings}
\label{sec_4_1_1}
\begin{figure}[t]
\centering
\includegraphics[scale=0.4]{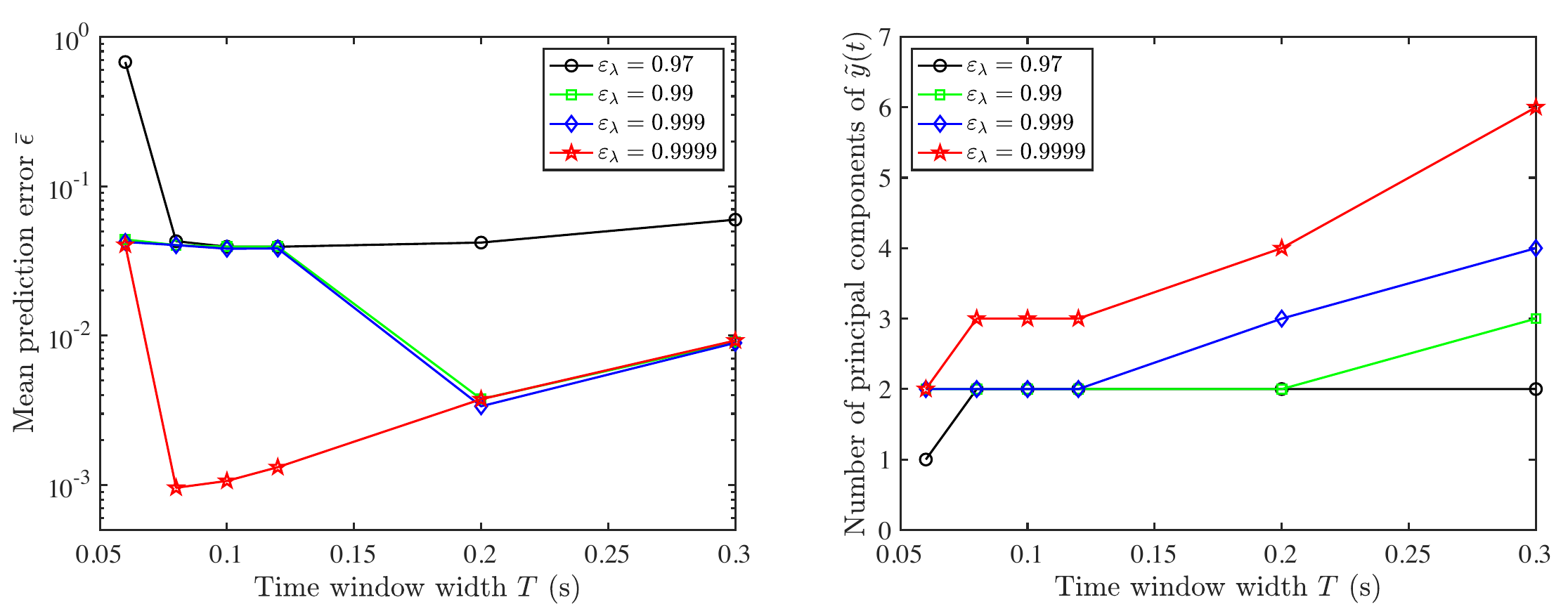}
\caption{Plots of prediction error and number of principal components of local response function $\tilde{y}(t)$ as functions of the time window width $T$ and retained variance proportion $\varepsilon_{\lambda}$.}
\label{Ex1_fig_diff_nW_varprop}
\end{figure}

Fig.~\ref{Ex1_fig_diff_nW_varprop} presents the 
mean prediction error and number of principal components of local response function $\tilde{y}(t)$ as functions of the time window width $T$ and retained variance proportion $\varepsilon_{\lambda}$. Here, the sizes of training dataset and testing dataset are 100 and 10,000 respectively. Both training and testing dataset are generated through MCS. As shown in the left panel of Fig.~\ref{Ex1_fig_diff_nW_varprop}, a general trend is observed where increasing the retained variance proportion $\varepsilon_{\lambda}$ results in lower prediction error. This is expected, as a higher $\varepsilon_{\lambda}$ retains more information for accurately predicting the dynamic response. Additionally, both excessively small and large values of the time window width $T$ lead to increased prediction errors. A too small window may fail to capture sufficient temporal information for model learning, whereas a too large window increases the complexity of the model, making it more difficult to learn effectively. The optimal parameters for F2NARX in this example are found to be $[T,\varepsilon_{\lambda}]=[0.08\ \text{s},~0.9999]$, and these settings will be used in the subsequent tasks for this example. The right panel of Fig.~\ref{Ex1_fig_diff_nW_varprop} shows that the number of principal components required to represent local response functions increases with both $\varepsilon_{\lambda}$ or $T$, which is intuitive, as more information needs to be retained with larger retained variance or longer window duration.

\begin{figure}[t]
\centering
\includegraphics[scale=0.4]{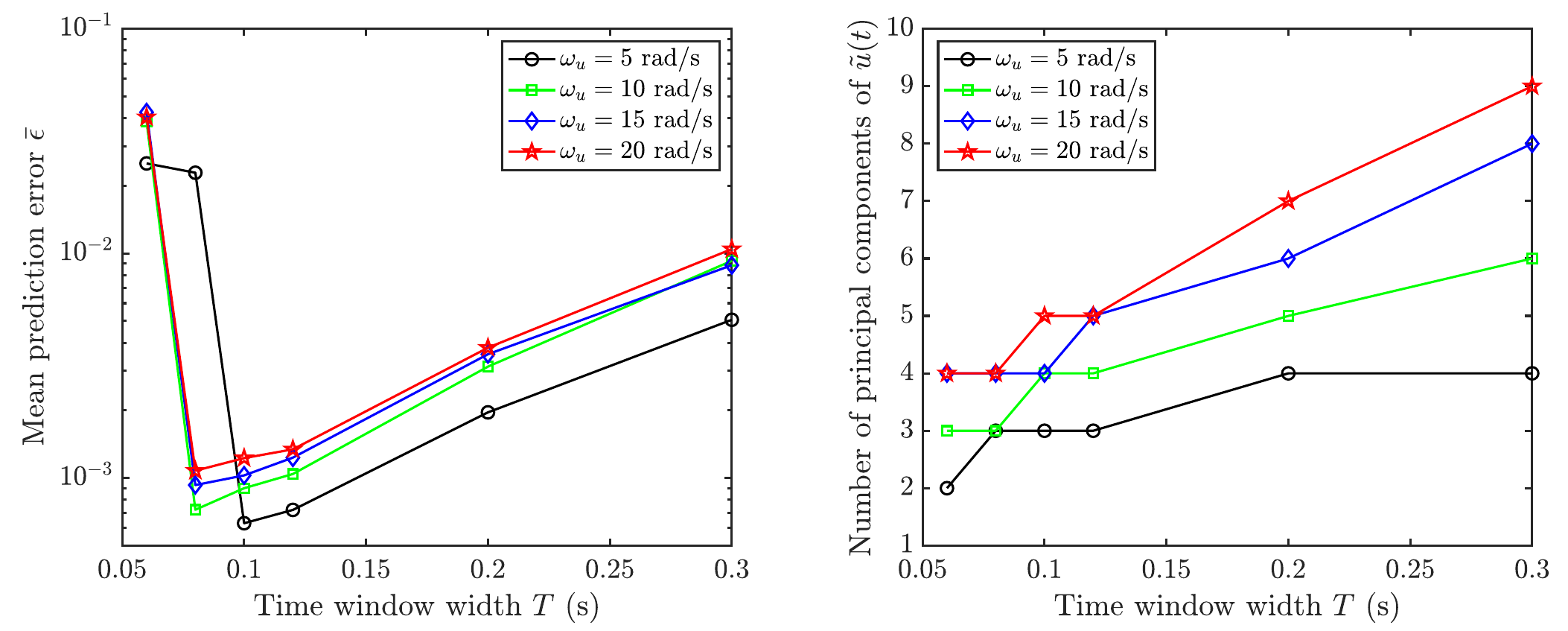}
\caption{Plots of prediction error and number of principal components of local excitation function $\tilde{u}(t)$ as functions of the time window width $T$ and upper cut-off frequency $\omega_u$.}
\label{Ex1_fig_diff_omega_T}
\end{figure}

To further investigate the influence of excitation frequency content, we fix $\varepsilon_{\lambda}$ at 0.9999 and evaluate the performance of F2NARX under different time window width $T$ and the excitation upper cut-off frequency $\omega_u$.
Fig.~\ref{Ex1_fig_diff_omega_T} presents the 
mean prediction error and number of principal components of local excitation function $\tilde{u}(t)$ as functions of $T$ and $\omega_u$.
As shown in Fig.~\ref{Ex1_fig_diff_omega_T}, when $\omega_u$ increases from 5 rad/s to 20 rad/s, the lowest prediction error shows a clear increasing trend.
This is because a larger $\omega_u$ introduces higher-frequency components into the excitation, and more principal components are required to characterize the local excitation function $\tilde{u}(t)$, as observed in the right panel of Fig.~\ref{Ex1_fig_diff_omega_T}.
This increases the dimensionality of the F2NARX model and makes the model more difficult to learn.
Another observation is that when $\omega_u$ is small, such as 5 rad/s, the time window width $T$ corresponding to the lowest prediction error is 0.1 s, which is larger than that for larger values of $\omega_u$.
This can be explained by the fact that excitations dominated by low-frequency components require a larger time window width to sufficiently capture their temporal patterns.

\subsubsection{Comparison of F2NARX with SGP-F-NARX}
\begin{figure}[t]
\centering
\includegraphics[scale=0.4]{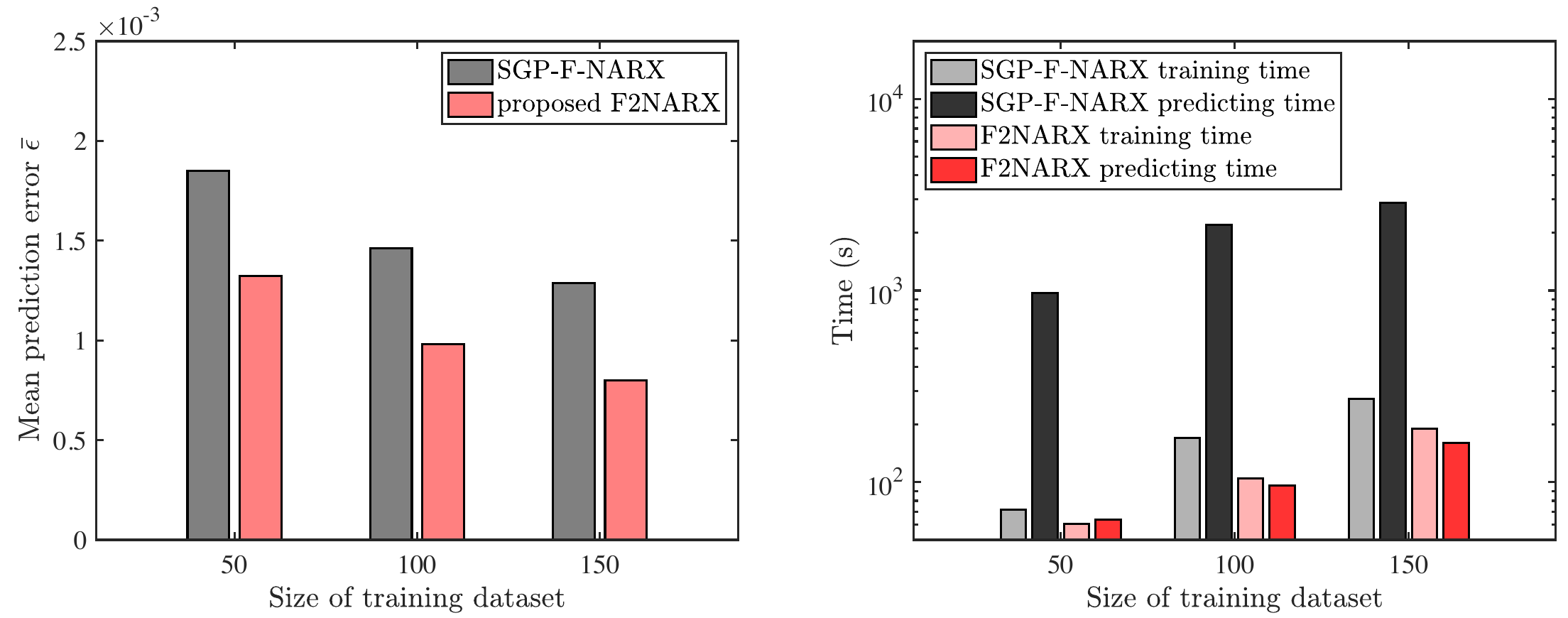}
\caption{Prediction error and computational time of SGP-F-NARX and the proposed F2NARX for the Bouc-Wen oscillator under different sizes of training dataset.}
\label{Ex1_fig_diff_methods}
\end{figure}

For the SGP-F-NARX method, a 5-fold cross-validation was performed on the training data set containing 50 training time histories.
The 50 time histories are divided into five subsets, each containing 10 time histories.
In each fold, one subset is used as the validation set, while the remaining four subsets are used to train the SGP-F-NARX model. 
This procedure is repeated five times so that each subset is used once for validation. 
For each candidate pair $[T,\varepsilon_{\lambda}]$, the prediction error is averaged over the five validation folds, and the parameter pair with the lowest average validation error is selected. 
The optimal parameters are found to be $[T,\varepsilon_{\lambda}]=[0.4\ \text{s},~0.99]$.
Fig.~\ref{Ex1_fig_diff_methods} presents the prediction error and computational time of both SGP-F-NARX and the proposed F2NARX method on the same testing dataset of size 10,000 under varying training dataset sizes. 

As shown in the left panel of Fig.~\ref{Ex1_fig_diff_methods}, the mean prediction error $\bar{\epsilon}$ decreases for both methods as the training dataset size increases. Across all training set sizes, the F2NARX method consistently yields lower prediction errors than SGP-F-NARX, although the prediction errors are generally within the same order of magnitude. The right panel of Fig.~\ref{Ex1_fig_diff_methods} highlights a significant advantage of F2NARX in prediction time, demonstrating more than one order of magnitude reduction compared to SGP-F-NARX. This improvement is attributed to the window-by-window prediction mechanism employed in the F2NARX framework.

\subsubsection{Results of probabilistic prediction} 
\begin{figure}[t]
\centering
\includegraphics[scale=0.4]{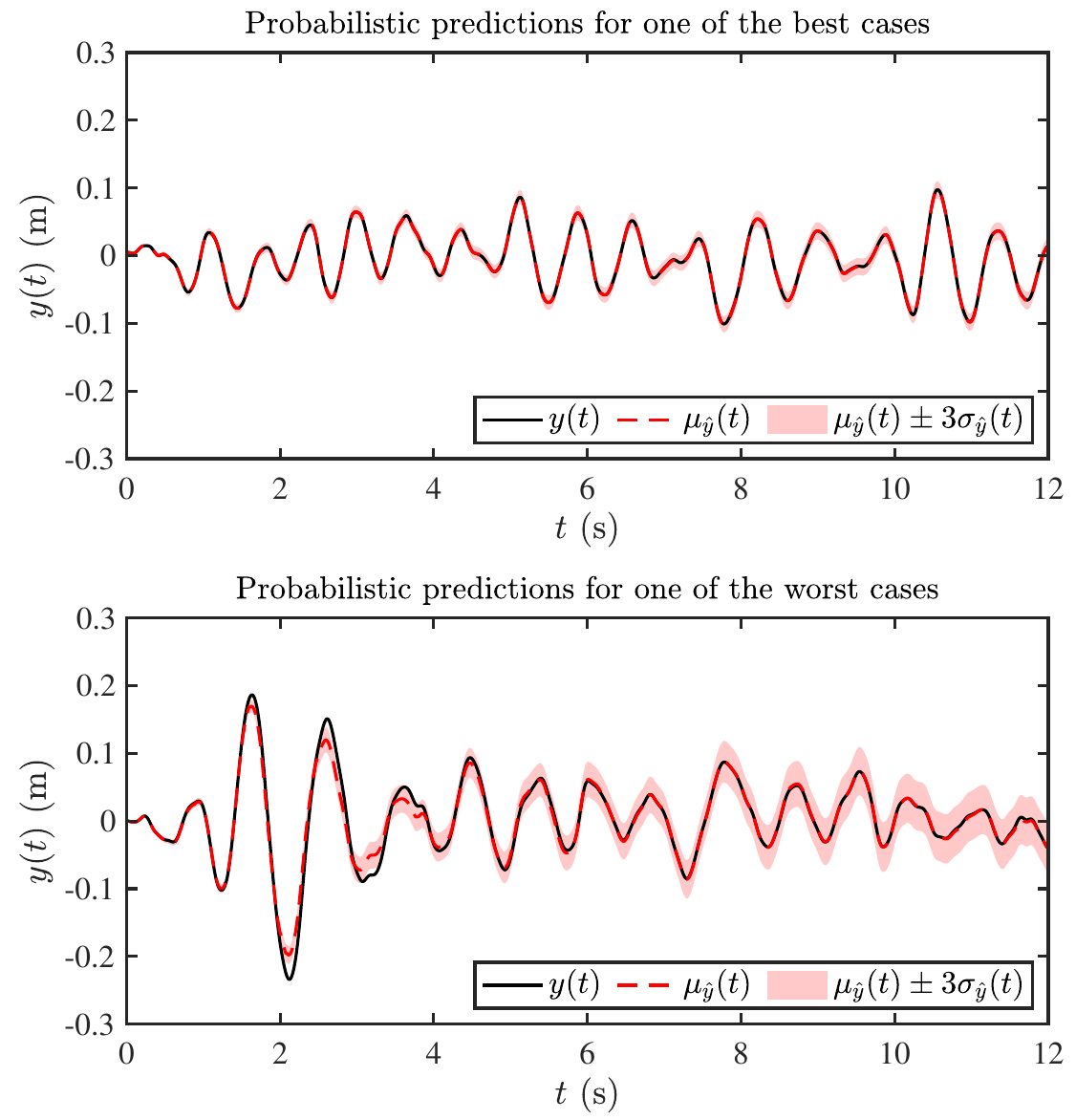}
\caption{Probabilistic predictions for one of the best cases and the worst cases for the Bouc-Wen oscillator.}
\label{Ex1_prob_prediction}
\end{figure}

Table~\ref{Ex1_prob_results} reports the mean prediction error of the prediction standard deviation function and the corresponding computational time for different probabilistic prediction methods, using an F2NARX model trained with 50 time trajectories and evaluated on 10,000 test trajectories.
The two approaches yield similarly small prediction errors. In terms of computational time, the proposed unscented transform-based method requires approximately half the time of the Taylor expansion-based method, and both approaches achieve nearly a three-order-of-magnitude reduction in computational cost compared to MCS. These results highlight the accuracy and efficiency of the proposed probabilistic prediction method. Furthermore, Fig.~\ref{Ex1_prob_prediction} shows the probabilistic predictions generated by F2NARX for cases with both small and large prediction errors.

\begin{table}[htbp]
  \centering
  \caption{Prediction error and computational time for different probabilistic prediction methods.}
  \label{Ex1_prob_results}
  \begin{tabular}{l l l}
    \toprule
    Methods & Prediction error  & Computational time (s) \\
    \midrule
    MCS & -- & $1.38\times10^5$  \\
    Taylor Expansion &  $8.39\times10^{-3}$  & $4.45\times10^2$ \\
    Unscented Transform &  $8.65\times10^{-3}$ & $2.58\times10^2$  \\
    \bottomrule
  \end{tabular}
\end{table}

\subsubsection{Results of active learning}
\begin{figure}[t]
\centering
\includegraphics[scale=0.4]{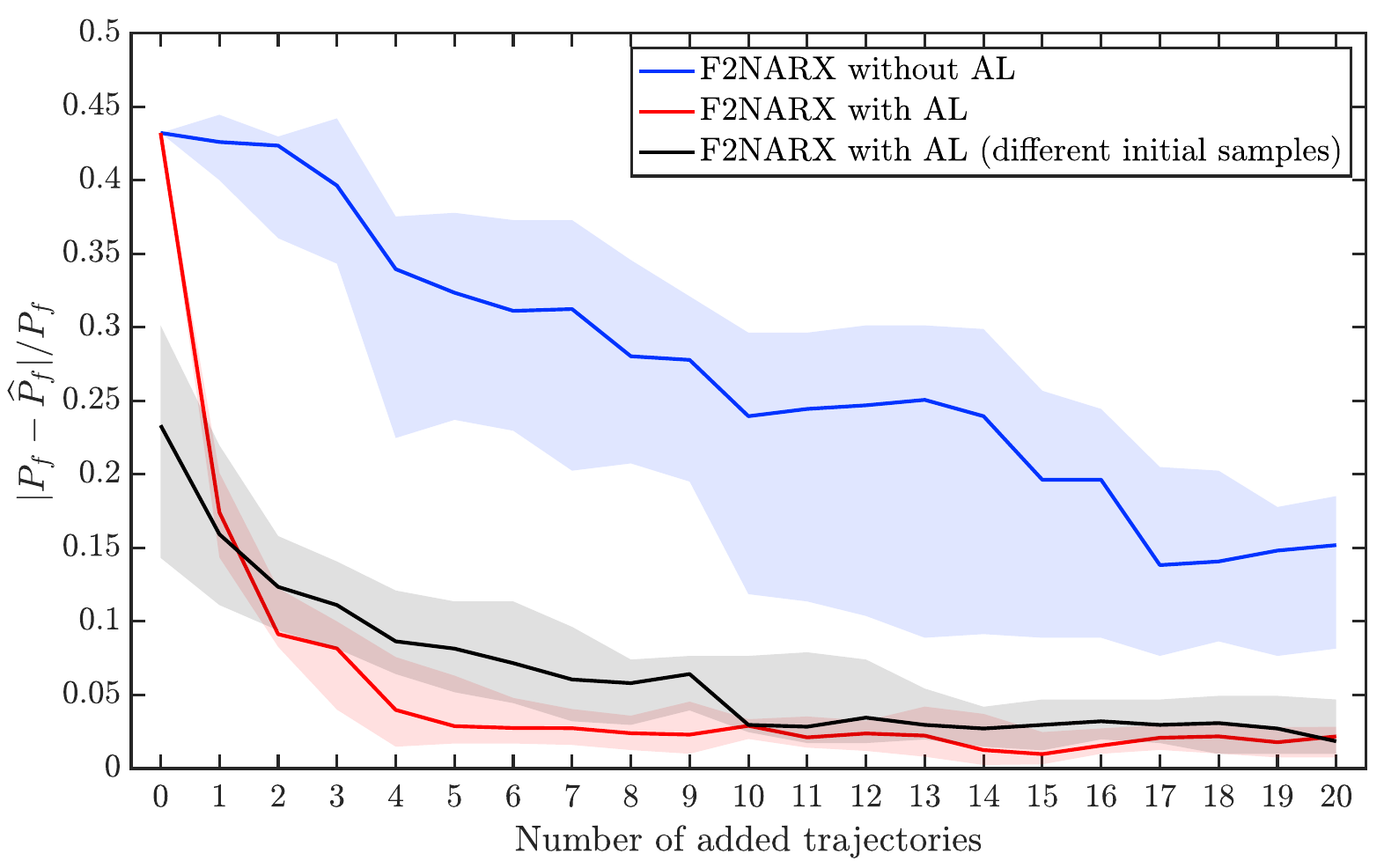}
\caption{Failure probability estimation error vs. iterations for the Bouc–Wen oscillator, with and without active learning.}
\label{Ex1_AL}
\end{figure}

In this example, failure is defined to occur when the maximum value of $|y(t)|$ exceeds $0.14\ \text{m}$. 
Accordingly, the dynamic performance function is given by $g(t)=0.14-|y(t)|$. 
The true first-passage failure probability is around 0.0405, the size of the MCS sample pool for active learning is set to $10^4$. 
Fig.~\ref{Ex1_AL} shows the evolution of the first-passage failure probability estimation error with and without active learning, labeled as F2NARX with AL and F2NARX without AL, respectively.
Both approaches start with the same initial training dataset of 10 trajectories, followed by 20 additional trajectories added either through the learning function or by random selection from the MCS sample pool. 
Each approach is repeated 10 times with the same initial training dataset but different MCS sample pools across replicates.
Within each replicate, the two approaches share the same MCS sample pool and are validated on the same MCS dataset.
To investigate the sensitivity to the initial training trajectories, the active learning approach is also repeated 10 times using different initial training datasets but the same MCS sample pool, labeled as F2NARX with AL (different initial samples).
The bold lines denote the median estimation errors, while the lower and upper boundaries of the shaded bands represent the 25th and 75th percentiles, respectively.

As shown in Fig.~\ref{Ex1_AL}, when additional trajectories are added, the active learning strategy rapidly reduces the estimation error, with the median curve converging after only five iterations.
In contrast, the non-active learning approach exhibits a much slower reduction in error and considerably larger variability across replicates. 
Furthermore, although larger variability is observed when active learning starts from different initial training trajectories, its advantage over the approach without active learning remains evident.
These results demonstrate that the probabilistic F2NARX predictions provide useful uncertainty information for guiding active learning in the considered first-passage probability estimation problem.

\subsection{Example 2: Three-story steel frame under seismic loading}
\begin{figure}[t]
\centering
\includegraphics[scale=0.6]{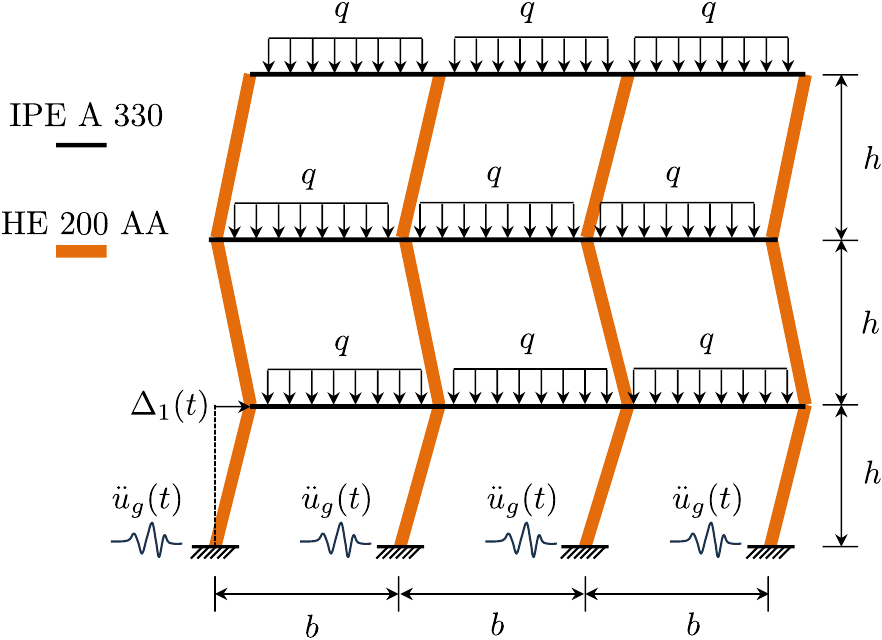}
\caption{Illustration of the three-story steel frame under seismic loading (figure adopted from \cite{schar2025surrogate}).}
\label{Ex2_steel_frame}
\end{figure}

The second example considers a realistic application of a nonlinear three-story steel frame under seismic loading \cite{zhu2023seismic,schar2025surrogate}. Fig.~\ref{Ex2_steel_frame} presents an illustration of the frame structure. We refer to \cite{zhu2023seismic} for more details of the model and \cite{schar2025surrogate} for the geometric parameters, material properties, and the live load of the structure. The quantity of interest is the interstory drift of the first floor $\Delta_1(t)$. 

The exogenous input is the ground acceleration $\ddot{u}_g(t)$, which is modeled by spectral representation method as \cite{deodatis2025spectral}:
\begin{equation}
    \label{ground_acc}
    \ddot{u}_{g}(t)=\sqrt{2}\sum_{i=1}^N \sqrt{2A^2(\omega_i,t)S_{\ddot{u}}(\omega_i)\Delta\omega}\cos (\omega_i t+\phi_i),
\end{equation}
where $\Delta\omega=\omega_u/N$ and $\omega_u$ is the upper cut-off frequency; $\omega_i=i\Delta\omega$; $\phi_1,\phi_2,\cdots,\phi_N$ are independent uniform random variables in $[0,2\pi]$ and $N$ is set to 1,000; $S_{\ddot{u}}(\omega)$ is the power spectral density of the ground motion, which is modeled by the Clough-Penzien spectrum as \cite{clough1975dynamics}:
\begin{equation}
    \label{CP}
    S_{\ddot u}(\omega)=
    \frac{\omega_g^4 + 4\zeta_g^2 \,\omega_g^2\,\omega^2}
     {(\omega_g^2 - \omega^2)^2 + \bigl(2\zeta_g\,\omega_g\,\omega\bigr)^2}\;
     \frac{\omega^4}
     {(\omega_f^2 - \omega^2)^2 + \bigl(2\zeta_f\,\omega_f\,\omega\bigr)^2}
     S_0,
\end{equation}
in which $\omega_g$ and $\zeta_g$ are the frequency and damping ratio of site soil, respectively; $\omega_f$ and $\zeta_f$ are the frequency and damping ratio of the high-pass filter, respectively; $S_0$ is the intensity factor of ground motions. $A(\omega,t)$ in Eq.~(\ref{ground_acc}) is the time-frequency modulation function:
\begin{equation}
    \label{modulation_function}
    A(\omega, t)=\exp\left(-\eta_0\frac{\omega t}{\omega_{g}T_g}\right)
    \left[\frac{t}{c}\exp\left(1-\frac{t}{c}\right)\right]^{r},
\end{equation}
where $\eta_0$ is the frequency modulation factor; $T_g$ is the time duration of ground motion; $c$ is the approximate arrival time of peak ground acceleration; $r$ is the shape control coefficient. The parameters in the spectral representation are set as follows: $\omega_u=50\pi~ \text{rad/s}$, $\omega_g=15~ \text{rad/s}$, $\omega_f=1.5~ \text{rad/s}$, $\zeta_g=\zeta_f=0.6$, $\eta_0=0.15$, $r=2$. To induce greater variability in ground acceleration, $S_0$ is assumed to follow a lognormal distribution with $\mu_{S_0}=\sigma_{S_0}=1.5\times10^{-6}~\text{m}^2/\text{s}^3$, and $c$ is modeled as a uniformly distributed random variable in the interval $[1,15]~\text{s}$. The time duration of the ground motion is set to $T_g=45~\text{s}$, which is discretized into 9,001 equally spaced time instants with time increment $\delta t=0.005~\text{s}$.

\begin{figure}[t]
\centering
\includegraphics[scale=0.4]{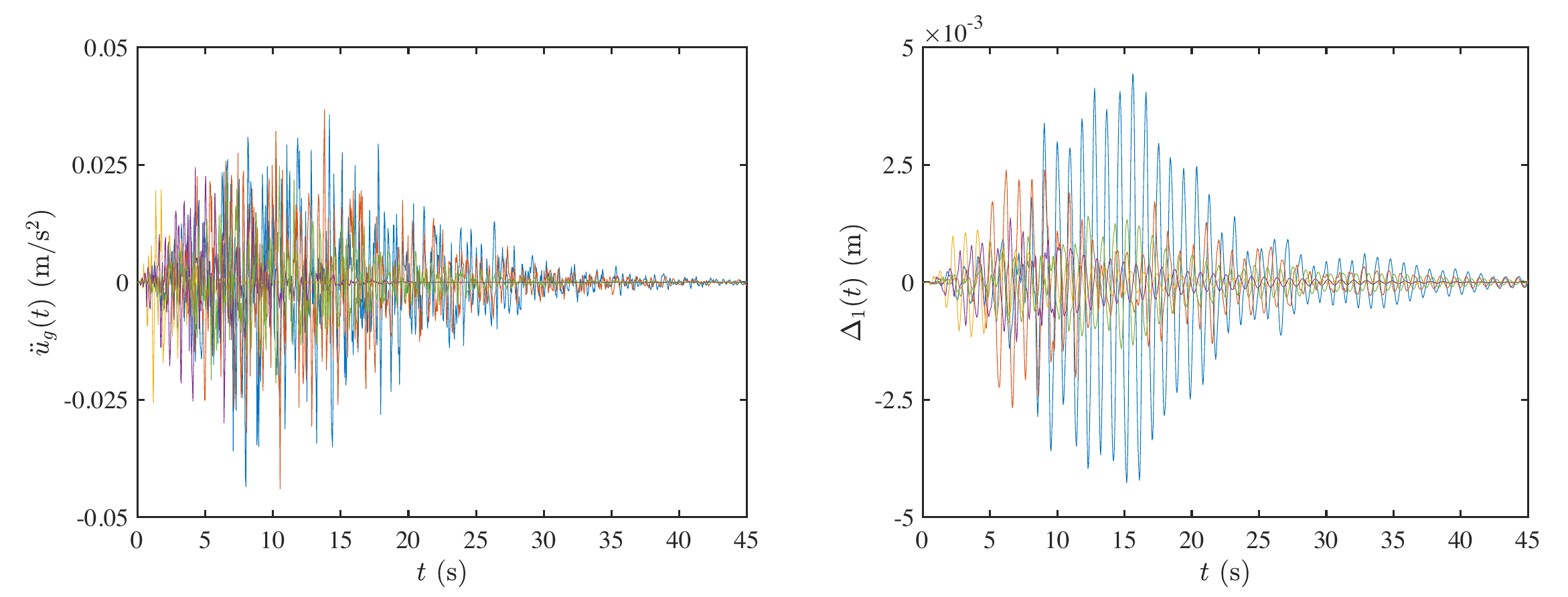}
\caption{Five different realizations of excitations $\ddot{u}_g(t)$ and corresponding responses $\Delta_1(t)$ for the three-story steel frame.}
\label{Ex2_u(t)_y(t)}
\end{figure}

A total of 10,000 ground motions were generated using MCS, and the corresponding interstory drift time histories $\Delta_1(t)$ were obtained with the open-source finite element software OpenSees \cite{mckenna2010nonlinear}. Running these 10,000 simulations in OpenSees required approximately 30.6 hours. Fig.~\ref{Ex2_u(t)_y(t)} shows five different realizations of ground motion and first interstory drift for this example.

\subsubsection{Performance of F2NARX method under different parameter settings}
\begin{figure}[t]
\centering
\includegraphics[scale=0.4]{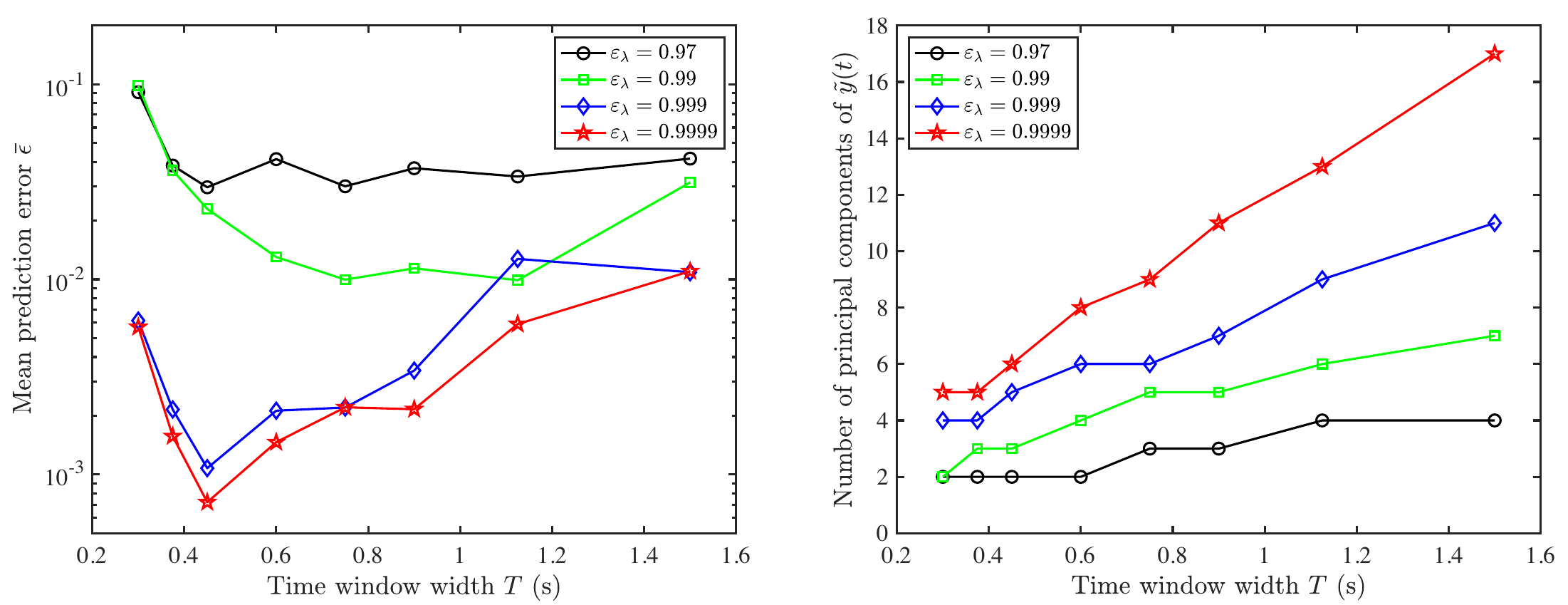}
\caption{Plots of prediction error and number of principal components of local response function $\tilde{y}(t)$ as functions of the time window width $T$ and retained variance proportion $\varepsilon_{\lambda}$ for the three-story steel frame.}
\label{Ex2_fig_diff_nW_varprop}
\end{figure}

Fig.~\ref{Ex2_fig_diff_nW_varprop} presents the prediction error and number of principal components of local response function $\tilde{y}(t)$ as functions of the time window width $T$ and retained variance proportion $\varepsilon_{\lambda}$. Here, the 10,000 trajectories obtained from the FEM simulations were randomly divided into a training dataset consisting of 50 trajectories and an out-of-sample test set of 9950 trajectories for performance evaluation. As shown in the left panel of Fig.~\ref{Ex2_fig_diff_nW_varprop}, a general trend is again observed where increasing the retained variance proportion $\varepsilon_{\lambda}$ results in lower prediction error. Additionally, both excessively small and large values of $T$ result in increased prediction errors. The optimal parameters for F2NARX in this example are found to be $[T,\varepsilon_{\lambda}]=[0.45\ \text{s},~0.9999]$, which are employed in the subsequent tasks. As shown in the right panel of Fig.~\ref{Ex2_fig_diff_nW_varprop}, the number of principal components required to represent the local response functions increases with either $\varepsilon_{\lambda}$ or $T$, since more information needs to be retained.

\subsubsection{Comparison of F2NARX with SGP-F-NARX}
\begin{figure}[t]
\centering
\includegraphics[scale=0.4]{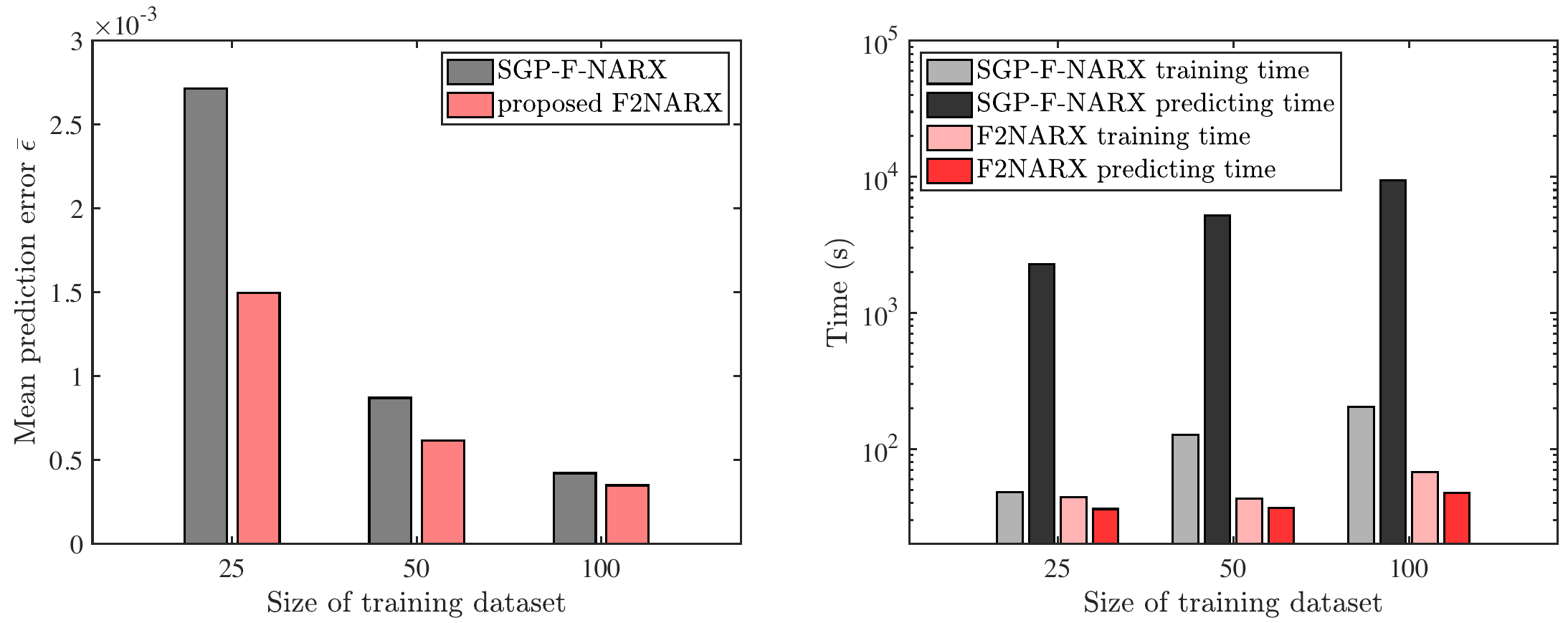}
\caption{Prediction error and computational time of SGP-F-NARX and the proposed F2NARX for the three-story steel frame under different sizes of training dataset.}
\label{Ex2_fig_diff_methods}
\end{figure}

For the SGP-F-NARX method, a 5-fold cross-validation was performed on the training data set containing 50 training time histories to select the optimal parameters, which are found to be $[T,\varepsilon_{\lambda}]=[1\ \text{s},~0.99]$.
Fig.~\ref{Ex2_fig_diff_methods} presents the prediction error and computational time of both SGP-F-NARX and the proposed F2NARX method under varying training dataset sizes. For each training dataset size, the training samples are randomly selected from the 10,000 trajectories, and the remaining trajectories are used for performance evaluation.

As shown in the left panel of Fig.~\ref{Ex2_fig_diff_methods}, the mean prediction error $\bar{\epsilon}$ decreases for both methods as the training dataset size increases. Across all training set sizes, both SGP-F-NARX and F2NARX produce prediction errors of the same order of magnitude, while F2NARX consistently achieves lower errors. The right panel of Fig.~\ref{Ex2_fig_diff_methods} highlights a significant advantage of F2NARX in prediction time, demonstrating about two orders of magnitude reduction compared to SGP-F-NARX. Specifically, F2NARX is 63, 142, and 198 times faster than SGP-F-NARX when making predictions with 25, 50, and 100 training samples, respectively. Furthermore, with 50 training samples, F2NARX is approximately 3,000 times faster than the finite element model.

\subsubsection{Results of probabilistic prediction}
\begin{figure}[t]
\centering
\includegraphics[scale=0.4]{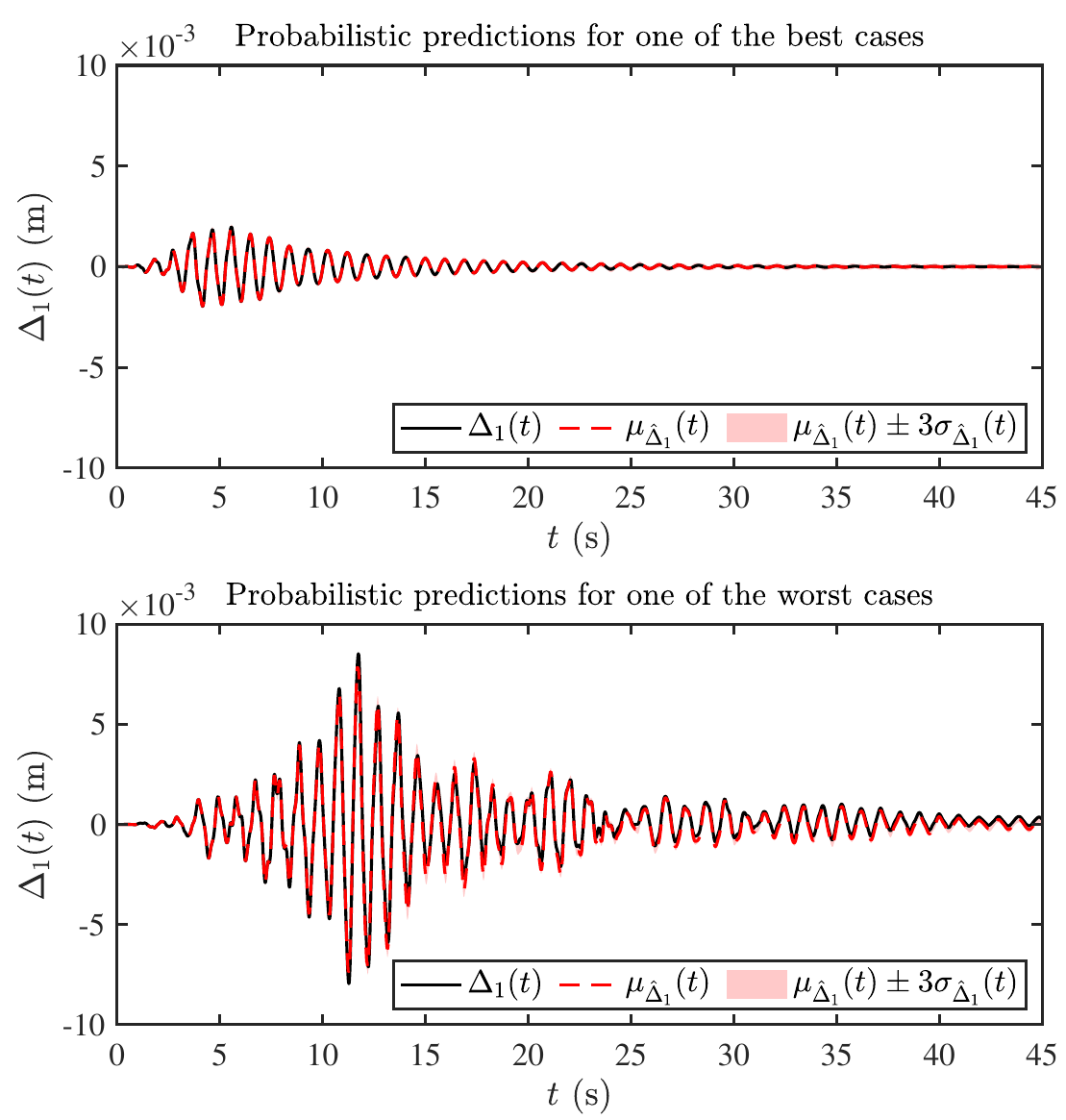}
\caption{Probabilistic predictions for one of the best cases and the worst cases for the three-story steel frame.}
\label{Ex2_prob_prediction}
\end{figure}

Table~\ref{Ex2_prob_results} summarizes the mean prediction errors of the predicted standard deviation function and the corresponding computational times for different probabilistic prediction methods. Both probabilistic prediction methods achieve similarly low prediction errors. However, the unscented transform-based approach demonstrates greater computational efficiency than the Taylor expansion-based approach. Moreover, both approaches offer nearly two orders of magnitude reduction in computational cost compared to MCS. These results highlight the accuracy and efficiency of the proposed probabilistic prediction strategy. Additionally, Fig.~\ref{Ex2_prob_prediction} illustrates the probabilistic predictions generated by F2NARX for cases with both small and large prediction errors.

\begin{table}[htbp]
  \centering
  \caption{Prediction error and computational time for different probabilistic prediction methods for the three-story steel frame.}
  \label{Ex2_prob_results}
  \begin{tabular}{l l l}
    \toprule
    Methods & Prediction error  & Computational time (s) \\
    \midrule
    MCS & -- & $2.13\times10^5$  \\
    Taylor Expansion &  $1.10\times10^{-3}$  & $2.09\times10^3$ \\
    Unscented Transform &  $1.09\times10^{-3}$ & $9.35\times10^2$  \\
    \bottomrule
  \end{tabular}
\end{table}

\subsubsection{Results of active learning}
\begin{figure}[t]
\centering
\includegraphics[scale=0.4]{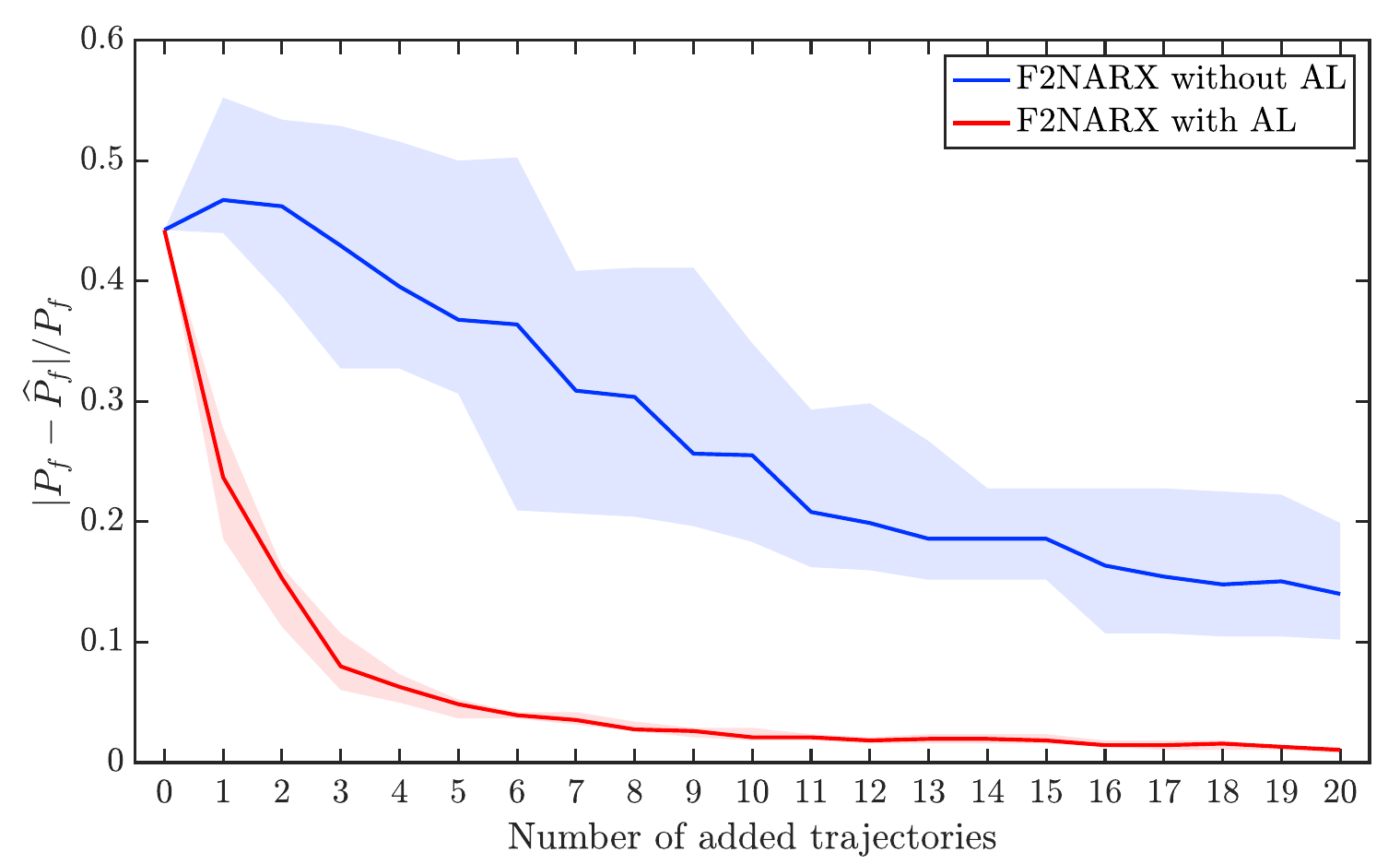}
\caption{Failure probability estimation error vs. iterations for the three-story steel frame, with and without active learning.}
\label{Ex2_AL}
\end{figure}

In this example, failure is defined to occur when the maximum value of $|\Delta_1(t)|$ exceeds $3.8\times10^{-3}\ \text{m}$. Accordingly, the dynamic performance function is given by $g(t)=3.8\times10^{-3}-|\Delta_1(t)|$. The true first-passage failure probability is around 0.0411, the size of the MCS sample pool for active learning is set to $10^4$. Fig.~\ref{Ex2_AL} presents the iteration curves of failure probability estimation error with and without active learning. Both methods begin with the same initial training dataset of 10 trajectories, but their behaviors significantly diverge as additional trajectories are added. The active learning strategy rapidly reduces the error, reaching values below 0.05 after only about five iterations, and then stabilizes with consistently low error and narrow variability across replicates. In contrast, the non--active learning approach shows a much slower decrease in error, with median values remaining above 0.1 after adding 20 trajectories, and with large variability across replicates. 
These results clearly demonstrate that the probabilistic F2NARX model provides useful uncertainty information for 
guiding active learning in the first-passage probability estimation problem.

\section{Conclusions and outlook}
\label{sec_5}
In this paper, we introduce a function-on-function nonlinear autoregressive model with exogenous inputs (F2NARX) as a surrogate modeling approach for complex dynamical systems. 
By extending the recently proposed $\mathcal{F}$-NARX formulation to a function-on-function autoregressive paradigm, F2NARX substantially improves prediction efficiency through a one-time-window-ahead strategy instead of the conventional one-step-ahead approach.
Through principal component analysis, the excitation and response functions within local time windows are represented by uncorrelated features, thereby decomposing the function-on-function mapping into a series of single-output mappings. 
Sparse Gaussian process regression is then employed to learn these mappings efficiently under large autoregressive training datasets. 
Furthermore, F2NARX leverages the unscented transform to enable efficient probabilistic prediction of dynamical responses based on the prediction uncertainty of SGP models. 
To further demonstrate the usefulness of these probabilistic predictions, F2NARX is combined with active learning for first-passage probability estimation.
A numerical example and a complex engineering case involving a finite element model are investigated to evaluate the performance of F2NARX. 
The results indicate that F2NARX can emulate dynamical responses with up to orders-of-magnitude reductions in computational time compared with state-of-the-art NARX models, while achieving higher accuracy. In the engineering case, F2NARX trained with only 50 time histories requires just 37 seconds to predict 10,000 dynamical responses of a nonlinear multi-degree-of-freedom steel frame structure on a conventional laptop, approximately three orders of magnitude faster than the finite element model and over two orders of magnitude faster than state-of-the-art NARX models, while maintaining a normalized mean squared error below $1\times10^{-3}$. 
For the moderate failure probabilities considered in this study, active learning guided by the probabilistic F2NARX predictions provides accurate and robust first-passage probability estimates with only about 15 training trajectories.

Two key parameters that significantly influence the performance of F2NARX are the time window width $T$ and the retained variance proportion $\varepsilon_{\lambda}$ for principal component analysis. Since the training and prediction process of the F2NARX model is highly efficient, cross-validation can be readily employed to determine suitable values of $T$ and $\varepsilon_{\lambda}$. Our results show that setting $\varepsilon_{\lambda}=0.9999$ generally works well. For the time window width, a good starting point is to vary $T$ within the range $[0.1\bar{T}_0,\ \bar{T}_0]$, where $\bar{T}_0$ denotes the average natural period corresponding to the lowest natural frequency of the system. In this study, we assume that the look-ahead time $T^+$ and look-back time $T^-$ are identical. However, for more complex dynamical systems, a longer look-back time and a shorter look-ahead time may be required. Future research could therefore investigate the use of different $T^+$ and $T^-$ values and develop methods to automatically determine them.
In addition, predictions in this work are made using non-overlapping time windows, which may cause discontinuities at the window boundaries. 
For problems where strict smoothness at window boundaries is critical, possible extensions include overlapping-window averaging or continuity-penalized training.

It should be noted that the proposed F2NARX model may also suffer from the curse of dimensionality if the external excitation or the response within local time windows exhibits large variability, requiring a large number of features for representation, if multiple external excitations are involved, or if the system parameters are high-dimensional. In such cases, advanced feature extraction techniques, such as autoencoders, may be integrated with F2NARX.
In addition, the current probabilistic prediction scheme relies on unscented transform. 
The accuracy of the UT may deteriorate when the underlying transformation is strongly nonlinear, the propagated distribution is significantly non-Gaussian, or the effective input dimension is high.
In such cases, sparse grid methods, cubature methods, or quasi-Monte Carlo can be used.
Moreover, the current work focuses on a single response QoI.
For multiple correlated QoIs, one F2NARX model can be constructed for each QoI, while the histories of all relevant QoIs can be included as autoregressive inputs to account for their dependence.
For large-scale QoIs, model order reduction can first be performed to obtain a low-dimensional latent representation, and the F2NARX model can then be constructed in the latent space.
The F2NARX model can also be constructed on a problem-specific exogenous input manifold \cite{schar2024emulating} to scale to applications with a large number of external excitations.

Although the proposed F2NARX provides an effective framework for probabilistic prediction, its application to active learning remains time-consuming. For instance, in the engineering case study, the active-learning procedure required approximately 4 hours to add 20 time trajectories. This computational burden becomes more pronounced for problems involving very small failure probabilities.
Future research may therefore explore more efficient simulation strategies beyond MCS for estimating small first-passage failure probabilities, or develop faster probabilistic prediction approaches.
In addition, future work may focus on developing more efficient active learning strategies, including learning functions tailored to dynamical systems and advanced stopping criteria such as budget-informed criterion \cite{zhang2025ice} that go beyond the current resource-based criterion.
Moreover, the current variance-based PCA truncation scheme may discard low-energy but important transient features, which could affect the accuracy of reliability analysis.
In future work, more advanced feature extraction methods, such as wavelet-based features, autoencoders, or other task-informed feature extraction techniques, may be considered to overcome the limitations of standard PCA.
Finally, aside from probabilistic surrogate modeling, F2NARX could potentially be extended to design optimization, control, and digital twins of complex dynamical systems.

\section*{CRediT authorship contribution statement}
\textbf{Zhouzhou Song}: Conceptualization, Methodology, Software, Visualization, Validation, Writing -- original draft, Funding acquisition. 
\textbf{Marcos A. Valdebenito}: Conceptualization, Writing -- review \& editing, Supervision. 
\textbf{Styfen Schär}: Software, Resources, Writing -- review \& editing.
\textbf{Stefano Marelli}: Resources, Writing -- review \& editing, Supervision.
\textbf{Bruno Sudret}: Writing -- review \& editing, Supervision.
\textbf{Matthias G.R. Faes}: Writing -- review \& editing, Supervision, Funding acquisition.

\section*{Declaration of Competing Interest}
The authors declare that they have no known competing financial interests or personal relationships that could have appeared to influence the work reported in this paper.

\section*{Acknowledgments}
Zhouzhou Song gratefully acknowledges the support of the Alexander von Humboldt Foundation for his postdoctoral fellowship. Matthias G.R. Faes gratefully acknowledges support from the Alexander von Humboldt Foundation through the Henriette Herz Scouting Program. Part of the research is funded by the State Key Laboratory of Disaster Reduction in Civil Engineering, Tongji University under project number SLDRCE24-02.


\bibliographystyle{elsarticle-num-names} 
\bibliography{reference}

\end{document}